\documentclass[12pt,twoside,reqno,openany]{amsart}
\usepackage{courier}
\usepackage{etoolbox}
\makeatletter
\patchcmd{\@settitle}{\uppercasenonmath\@title}{}{}{}
\patchcmd{\@setauthors}{\MakeUppercase}{}{}{}
\patchcmd{\section}{\scshape}{}{}{}
\makeatother
\usepackage{amsbsy,amscd,amsfonts,amsmath,amssymb,amsthm,color,
fancybox,fancyhdr,footmisc,graphics,graphicx,ifthen,mathrsfs,
multicol,pdfpages,rotating,times,wasysym}

\usepackage[dvipsnames,svgnames,x11names]{xcolor}

\usepackage[all]{xy}
\usepackage[utf8]{inputenc}
\usepackage[T1]{fontenc}
\usepackage{mathtools}
\newtagform{EngelLie}[\scriptstyle]{$}{$}
\makeatletter\newcommand{\leqnomode}{\tagsleft@true}
\newcommand{\reqnomode}{\tagsleft@false}\makeatother

\newtheorem{Theorem}[equation]{Theorem}
\newtheorem{Proposition}[equation]{Proposition}

\theoremstyle{definition}

\newtheorem{Definition}[equation]{Definition}

\newcommand{\aaux}{{\text{\usefont{T1}{qcs}{m}{sl}a}}}
\newcommand{\baux}{{\text{\usefont{T1}{qcs}{m}{sl}b}}}
\newcommand{\caux}{{\text{\usefont{T1}{qcs}{m}{sl}c}}}
\newcommand{\daux}{{\text{\usefont{T1}{qcs}{m}{sl}d}}}
\newcommand{\eaux}{{\text{\usefont{T1}{qcs}{m}{sl}e}}}
\newcommand{\faux}{{\text{\usefont{T1}{qcs}{m}{sl}f}}}
\newcommand{\gaux}{{\text{\usefont{T1}{qcs}{m}{sl}g}}}
\newcommand{\haux}{{\text{\usefont{T1}{qcs}{m}{sl}h}}}
\newcommand{\iaux}{{\text{\usefont{T1}{qcs}{m}{sl}i}}}
\newcommand{\jaux}{{\text{\usefont{T1}{qcs}{m}{sl}j}}}

\newcommand{\waux}{{\text{\usefont{T1}{qcs}{m}{sl}w}}}
\newcommand{\xaux}{{\text{\usefont{T1}{qcs}{m}{sl}x}}}
\newcommand{\yaux}{{\text{\usefont{T1}{qcs}{m}{sl}y}}}
\newcommand{\zaux}{{\text{\usefont{T1}{qcs}{m}{sl}z}}}

\definecolor{blue}{cmyk}{1.,1.,0.,0.63}
\definecolor{red}{cmyk}{0.,1.,1.,0.63}
\definecolor{green}{cmyk}{1.,0.,1.,0.63}
\definecolor{black}{cmyk}{1.,1.,1.,1.}

\makeatletter
\renewcommand{\@fnsymbol}[1]
{\ensuremath{\ifcase#1\or $*$\or $**$\or $***$\or $****$\or $*****$
\else\@ctrerr\fi}}
\makeatother

\numberwithin{equation}{section}

\newcommand{\Section}[1]{
\renewcommand{\thesection}{\bf\arabic{section}}
\section{#1}
\renewcommand{\thesection}{\arabic{section}}}

\newcommand{\style}[1]{\text{\footnotesize{\sf #1}}}

\renewcommand{\det}{\style{det}}

\renewcommand{\lim}{\style{lim}}

\renewcommand{\mod}{\style{mod}}

\renewcommand{\Re}{\style{Re}}

\newcommand{\isqrt}{{\scriptstyle{\sqrt{-1}}}}

\newcommand{\smallsum}[1]{
\underset{#1}{\raisebox{1pt}{$\sum$\,}}
}

\newcommand{\vf}{\vfill\end{document}}

\title[Holomorphic immersion]{Holomorphic immersions of bi-disks into 9 dimensional\\ real hypersurfaces with Levi signature $(2,2)$}
\author{Wei Guo {\sc Foo}
\qquad Jo\"el {\sc Merker}}
\begin{document}
\begin{abstract}
Inspired by an article of R. Bryant on holomorphic immersions of unit disks into Lorentzian CR manifolds, we discuss the application of Cartan's method to the question of the existence of bi-disk $\mathbb{D}^{2}$ in a smooth $9$-dimensional real analytic real hypersurface $M^{9}\subset\mathbb{C}^{5}$ with Levi signature $(2,2)$ passing through a fixed point. The result is that the lift to $M^{9}\times U(2)$ of the image of the bi-disk in $M^{9}$ must lie in the zero set of two complex-valued functions in $M^{9}\times U(2)$. We then provide an example where one of the functions does not identically vanish, thus obstructing holomorphic immersions. 
\end{abstract}
\maketitle
\markleft{}
\markright{}


\Section{\bf Introduction}
\label{introduction}

In the theory of exterior differential systems, \'{E}lie Cartan designed a constructive and powerful method to classify geometric objects under transformation, a process which he called the {\sl equivalence method} (or {\sl m\'{e}thode d'\'{e}quivalence} in French). There is an extensive literature that provides references to the equivalence method (see {Bryant-Chern-Gardner-Goldschmidt-Griffiths}  \cite{BCGGG-1991}, {Ivey-Landsberg}  \cite{Cartanforbeginners}, {Olver} \cite{Olver-1995}).

The application of the equivalence method to CR manifolds was studied by {\'{E}. Cartan }\cite{Cartan-1}, \cite{Cartan-2}, and later by {Tanaka} \cite{Tanaka} and {Chern-Moser} \cite{Chern-Moser-1974}. Much later, there are also other literature revolving around this subject ({Beloshapka-Ezhov-Schmalz} \cite{Beloshapka-Ezhov-Schmalz-2007}, {Ezhov-Isaev-Schmalz} \cite{Ezhov-Isaev-Schmalz-1999}, {Ezhov-McLaughlin-Schmalz} \cite{Ezhov-McLaughlin-Schmalz-2011},  {Isaev-Zaitsev} \cite{Isaev-Zaitsev-2013},   {Medori-Spiro} \cite{Medori-Spiro-2014}, {Schmalz-Spiro} \cite{Schmalz-Spiro-2006}, etc.), as well as recent interests in the effective aspects of the Cartan process (see articles by {Merker J.}, {Pocchiola S.}, and {Sabzevari M.}  \cite{Merker-Manuscript-2013-1}, \cite{Merker-Manuscript-2013-2}, \cite{Merker-Manuscript-2013-3}, \cite{Merker-Lie-2015}, \cite{Pocchiola-Manuscript-2014-1}, \cite{Pocchiola-Manuscript-2014-2}, \cite{Sabzevari-Merker-2014}, \cite{Sabzevari-Merker-Manuscript-2014}, \cite{Merker-Pocchiola-2018}) .

The objective of this paper is to study the existence of holomorphic varieties embedded into Levi non-degenerate real-analytic CR manifolds with different Levi signatures using this method. One early reference is \cite{Bryant-1982} where R. Bryant considered the Lorentzian case.

Let $M^{2n+1}$ be a real smooth manifold of dimension $2n+1$ with  $n\geqslant 1$. Let $T^{1,0}M$ be a sub-bundle of the complexified tangent space $\mathbb{C}TM=TM\otimes_{\mathbb{R}}\mathbb{C}$. Then the tuple $(M,T^{1,0}M)$ is called {\sl CR manifold} with the {\sl CR structure} $T^{1,0}M$ if the following conditions are satisfied:
\begin{enumerate}
\item the rank of the bundle $T^{1,0}M$ is $n$,
\item $T^{1,0}M\cap T^{0,1}M=\{0\}$ where $T^{0,1}M$ is the complex conjugate of $T^{1,0}M$, and 
\item the CR structure is integrable in the sense of Frobenius, meaning that $[T^{1,0}M,T^{1,0}M]\subseteq T^{1,0}M$.
\end{enumerate}

Throughout, the CR manifolds are assumed to be real analytic. As the problem in this paper is local in nature, the following theorem of Andreotti explains the advantage of this added assumption:

\begin{Theorem}[Andreotti, See \cite{Chen-Shaw}, Theorem 12.1.3]
Any real analytic CR manifold $(M,T^{1,0}M)$ of dimension $2n+1$ with $n\geqslant 1$ can be locally embedded as a real analytic  hypersurface in $\mathbb{C}^{n+1}$.
\end{Theorem}

Henceforth, let $(z_{1},\cdots,z_{n+1})$ with $z_{i}=x_{i}+\isqrt y_{i}$ be holomorphic coordinates, and let $J:T\mathbb{C}^{n+1}\rightarrow T\mathbb{C}^{n+1}$ be the standard complex structure 
\[
J(\partial_{x_{i}})=\partial_{y_{i}},\qquad 
J(\partial_{y_{i}})=-\partial_{x_{i}}.
\]
The real analytic CR manifold may be (locally) defined by a smooth real analytic real-valued function $M=\{\rho=0\}$, with the CR structure naturally inherited from the ambient  space 
\[
T^{1,0}M=\mathbb{C}TM\cap T^{1,0}\mathbb{C}^{n+1},
\]
where $\mathbb{C}TM:=TM\otimes_{\mathbb{R}}\mathbb{C}$. The local differential $1$-form 
$\isqrt\overline{\partial}\rho$ 
on $M$ is complex-valued and it has restriction  
$\theta:=\isqrt\overline{\partial}\rho|_{M}$ to $M$ {\em real-valued},
cf. \cite[page 25]{Jacobowitz-1990}. The extension of $\theta$ to $\mathbb{C}TM$ satisfies
\[
\ker\theta=T^{1,0}M\oplus T^{0,1}M.
\] 
This $1$-form can be used to define a Hermitian product on $T^{1,0}M$ at each $p\in M$ as follows: for any two vectors $X_{p}$ and $Y_{p}$ in $T_{p}^{1,0}M$, and for any two $T^{1,0}M$ sections $X$ and $Y$ with $X(p)=X_{p}$ and $Y(p)=Y_{p}$, let 
\begin{equation}
\begin{aligned}
{\sf Lev}_{p}:T_{p}^{1,0}M\times T_{p}^{1,0}M &\longrightarrow \mathbb{C}\\
(X_{p},Y_{p}) &\longmapsto -\theta([X,\bar{Y}])(p).
\end{aligned}
\end{equation}
This bilinear map is well-defined, independent of the choice of $X$ and $Y$ which respectively give $X(p)$ and $Y(p)$ at $p$, and is known as the {\sl Levi form}. In view of the Cartan-Lie formula, which states that 
\[
d\theta(X\wedge \bar{Y})=
X\theta(\bar{Y})-\bar{Y}\theta(X)-\theta([X,\bar{Y}])
=
-\theta([X,\bar{Y}]),
\]
the $2$-form $d\theta$ modulo $\theta$ is also known as the Levi form. With respect to a $T^{1,0}M$ frame $\{\mathcal{L}_{1},\cdots,\mathcal{L}_{n}\}$, the {\sl Levi matrix} at $p$ is a matrix with entries $l_{ij}={\sf Lev}_{p}(\mathcal{L}_{i}(p),\overline{\mathcal{L}_{j}(p)})$. Then the manifold $M$ is {\sl Levi non-degenerate} if the Levi matrix is non-singular everywhere.

\subsection{Adapted coframes, holomorphic immersion, and Sommer's theorem} Assuming that $M$ is a real analytic, Levi non-degenerate CR manifold, there is a $T^{1,0*}M$ co-frame $\{\alpha^{1},\dots,\alpha^{n}\}$ which diagonalises $d\theta$. By renumbering the indices if necessary, this $2$-form may be written as 
\[
d\theta
\equiv 
\isqrt
\big(\alpha^{1}\wedge\bar{\alpha}^{1}
+
\cdots
+
\alpha^{p}\wedge\bar{\alpha}^{p}
-
\alpha^{p+1}\wedge\bar{\alpha}^{p+1}
-
\cdots
-
\alpha^{n}\wedge\bar{\alpha}^{n}\big)\qquad 
\mod\ \theta.
\]
Therefore, the signature of the Levi form is $(n_{+},n_{-})=(p,n-p)$. Multiplying by $-1$ if necessary, it may be assumed that $n_{+}:= p$ is less than $n_{-}:=n-p$.  If $\varphi:\mathbb{D}^{k}\rightarrow M$ is a holomorphic immersion, the signature of the Levi form does not allow $k$ to be greater than $p$. For simplicity, let $p=k$, and let $(s_{1},\dots,s_{p})\in\mathbb{D}^{p}$ be holomorphic coordinates.

One of the steps in Cartan's equivalence method is prolongation. To illustrate this step, it is required to study the pushforward of the sections of the $T^{1,0}\mathbb{D}^{p}$ bundle. For each $i$, let $\{\mathcal{A}_{i}:\ 1\leqslant i\leqslant n\}$ be a change of $T^{1,0}M$-frame that is dual to $\alpha^{i}$, obtained from the Gram-Schmidt process. Since $\varphi$ is a holomorphic immersion, there exist $p$ linearly independent vector fields on $\varphi(\mathbb{D}^{p})$ which can be expressed as 
\[
\varphi_{*}\partial_{s_{i}}=\smallsum{1\leqslant j\leqslant n}f_{i,j}(s_{1},\dots,s_{p},\bar{s}_{1},\dots,\bar{s}_{p})\mathcal{A}_{j}\big|_{\varphi(\mathbb{D}^{p})}
\]
for certain real-analytic functions $f_{i,j}$, satisfying the vanishing condition
\begin{equation}\label{intro-sommer-1}
d\theta(\varphi_{*}\partial_{s_i}\wedge \varphi_{*}\partial_{\bar{s}_j})\equiv 0.
\end{equation}
This implies that over $\varphi(\mathbb{D}^{p})$, the distribution of vector spaces spanned by $\varphi_{*}\partial_{s_{i}}$ lies in the isotropic cone of the Levi form. In other words, for each $v\in \varphi_{*}T^{1,0}\mathbb{D}^{p}$,
\[
d\theta(v\wedge\bar{v})\equiv 0.
\]
The fact that the isotropic cone contains a distribution of $p$-dimensional vector space is a clue to the first prolongation process. From equation \eqref{intro-sommer-1}, a direct substitution results in
\[
\smallsum{1\leqslant j\leqslant p}
|f_{i,j}|^{2}
-
\smallsum{1\leqslant j\leqslant n-p}
|f_{i,p+j}|^{2}
=
0.
\]
At this stage, define
\[
U_{n_{-},n_{+}}
:=
\big\{A\in M_{n_{-}\times n_{+}}(\mathbb{C}):\ A^{*}A={\sf Id}_{n_{+}\times n_{+}}\big\}.
\]
A theorem of Sommer (see theorem \ref{Sommer's-theorem}) shows that the positive and the negative part are related by a matrix in $U_{n_{-},n_{+}}$,
{\footnotesize
\[
\left(
\begin{matrix}
f_{1,p+1} & \cdots & f_{p,p+1}\\
\vdots & \ddots & \vdots \\
f_{1,n} & \cdots & f_{p,n}
\end{matrix}
\right)
=
\underbrace{
\left(
\begin{matrix}
f_{1,p+1} & \cdots & f_{p,p+1}\\
\vdots & \ddots & \vdots \\
f_{1,n} & \cdots & f_{p,n}
\end{matrix}
\right)
\left(
\begin{matrix}
f_{1,1}  & \cdots & f_{p,1}\\
\vdots & \ddots & \vdots\\
f_{1,p} & \cdots & f_{p,p}
\end{matrix}
\right)^{-1}
}_{:={\sf U}\in U_{n_{-},n_{+}}}
\cdot
\left(
\begin{matrix}
f_{1,1}  & \cdots & f_{p,1}\\
\vdots & \ddots & \vdots\\
f_{1,p} & \cdots & f_{p,p}
\end{matrix}
\right).
\]
}
The fact that ${\sf U}^{*}{\sf U}$ is the identity matrix is due to $\xaux^{*}\cdot\xaux-({\sf U}\xaux)^{*}\cdot {\sf U}\xaux=0$ for any $\xaux\in\mathbb{C}^{p}$ since the isotropic cone contains a $p$-dimensional vector space.  If $n_{-}=n_{+}:=p$, then $U_{n_{-},n_{+}}=U(p)$ is the usual set of unitary matrices of size $p$.

The first prolongation process treats ${\sf U}$ as any matrix in $U_{n_{-},n_{+}}$ satisfying the {\sl lifting condition}. More precisely, for any holomorphic immersion $\varphi:\mathbb{D}^{p}\rightarrow M$  into a CR real hypersurface with Levi signature $(p,n-p)$, there is a lift $\tilde{\varphi}:\mathbb{D}^{p}\rightarrow M\times U_{n_{-},n_{+}}$ that sends every point $p\in\mathbb{D}^{p}$ to $(\varphi(p),{\sf U})$ so that the following diagram commutes
\[
\xymatrix{
 & M\times U_{n_{-},n_{+}}\ar^{\pi}[d]\\
 \mathbb{D}^{p}\ar^{\tilde{\varphi}}[ur]\ar_{\varphi}[r]  & M.
}
\]
The map $\pi$ is just the projection onto the first component $M$. 

\subsection{The Pfaffian system}
Let ${\sf u_{i,j}}$ denote the coefficient of the $U_{n_{-},n_{+}}$ matrix ${\sf U}$. Consider the Pfaffian system, which is a system of differential $1$-forms:
\begin{equation}
\begin{aligned}
\omega^{0} &:= \theta,\\
\omega^{1} &:= \alpha^{1},\\
&\vdots\\
\omega^{p} &:= \alpha^{p},\\
\omega^{p+1} &:= 
\alpha^{p+1}-\smallsum{1\leqslant k\leqslant p}{\sf u_{1,k}}\alpha^{k},\\
&\vdots \\
\omega^{n} &:= 
\alpha^{n}-\smallsum{1\leqslant k\leqslant p}{\sf u_{n-p,k}}\alpha^k.
\end{aligned}
\end{equation}
 These differential $1$-forms constitute a  $T^{1,0}M$ co-frame over $M$. Let $\mathcal{I}$ be the ideal generated by $\omega^{0}$,  $\omega^{k}$  and $\bar{\omega}^{k}$ for $p+1\leqslant k\leqslant n$. The ideal then describes the bundle $\varphi_{*}T^{1,0}\mathbb{D}^{p}$  over $\varphi(\mathbb{D}^{p})$, and has to be expanded to a larger ideal $\mathcal{I}_{+}$ consisting of some differential $1$-forms on $M\times U_{n_{-},n_{+}}$ so that $d\mathcal{I}_{+}\equiv 0\,\,\mod\mathcal{I}_{+}$. In the CR-Lorentzian case \cite{Bryant-1982}, this allows certain values of ${\sf u_{i,j}}$ to be computed, giving the possible directions of the tangent vectors to $\varphi(\mathbb{D}^{p})$.

\subsection{A summary of Bryant's approach \cite{Bryant-1982}}

We first give a brief explanation of \cite{Bryant-1982}, restricting to the real hypersurfaces in $\mathbb{C}^{3}$ whose signature of the Levi form $d\theta$ is $(1,1)$ everywhere. The $2$-form can simply be expressed in terms of the adapted coframe as:
\[
d\theta\equiv 
\isqrt\big(\alpha^{1}\wedge\overline{\alpha}^{1}-\alpha^{2}\wedge\overline{\alpha}^{2}\big)\qquad \mod\ \theta.
\]

He then considered the problem of existence of holomorphic immersion of the unit disk $\mathbb{D}$ into $M$. Let $\varphi:\mathbb{D}\rightarrow M^{5}$ be such an immersion. Denoting by $t$ a holomorphic coordinate of $\mathbb{D}$, the holomorphic tangent vector to the image $\varphi(\mathbb{D})$ may be written as 
\[
\mathcal{L}_{\varphi(t)}=f_{1}(t,\bar{t})\mathcal{A}_{1}+f_{2}(t,\bar{t})\mathcal{A}_{2}.
\]
This vector field lies in the isotropic cone of the Levi form, and hence
\[
|f_{1}|^{2}-|f_{2}|^{2}=0.
\]
The hypothesis  that  $\varphi$ is an immersion does not allow $f_{1}$ or $f_{2}$ to vanish at any point in $M^{5}$.  If $f_{1}$ is nowhere vanishing, then 
\[
\left|\frac{f_{2}}{f_{1}}\right|=1,
\]
and hence there exists a circle-valued function $\lambda: M^{5}\rightarrow \mathbb{S}^{1}$ such that $f_{2}=\lambda f_{1}$, inducing a lift $\tilde{\varphi}:\mathbb{D}\rightarrow M\times \mathbb{S}^{1}$ to the product space with fibre $\mathbb{S}^{1}$ so that $\pi\circ \tilde{\varphi}=\varphi$. Treating $\lambda$ as an unknown variable as part of the Cartan process, the following Pfaffian system is set up as before:
\begin{equation}
\begin{aligned}
\omega^{0} &:= \theta,\\
\omega^{1} &:= \alpha^{1},\\
\omega^{2} &:= \alpha^{2}-\lambda\alpha^{1}.
\end{aligned}
\end{equation}
It can be seen that $\varphi^{*}\theta\equiv 0\equiv \varphi^{*}\omega^{2}$ with $\varphi^{*}(\alpha^{1}\wedge\bar{\alpha}^{1})\neq 0$. Hence the ideal $\mathcal{I}$ generated by $\theta$, $\omega^{2}$ and $\bar{\omega}^{2}$ describes the $\varphi_{*}\mathbb{C}T\mathbb{D}$-bundle over $\varphi(\mathbb{D})$.  Since $\varphi^{*}$ commutes with the Poincar\'{e} differential operator $d$, it necessarily follows that $\varphi^{*}(d(\bar{\lambda}\omega^{2}))\equiv 0\,\mod\mathcal{I}$.  A direct calculation then shows that there exists a function $L$ on $M^{5}$ with
\[
d(\bar{\lambda}\omega^{2})=-(\bar{\lambda}d\lambda+L\bar{\alpha}^{1})\wedge\alpha^{1}\equiv 0\qquad \mod\ \mathcal{I}.
\]
Pulling back the equation to the disk by $\varphi$, 
\[
\varphi^{*}(\bar{\lambda}d\lambda+L\bar{\alpha}^{1})\wedge
\varphi^{*}\alpha^{1}\equiv 0,
\]
and so by Cartan's lemma, a modification of $\varphi^{*}(\bar{\lambda}d\lambda+L\bar{\alpha}^{1})$ by a multiple of $\varphi^{*}\alpha^{1}$  is required in order for the following $1$-form to vanish on $\varphi(\mathbb{D})$:
\[
\varphi^{*}(\bar{\lambda}d\lambda+L\bar{\alpha}^{1})+\mu\varphi^{*}\alpha^{1}\equiv 0.
\]
 When this happens, on $\varphi(\mathbb{D})$,
\begin{equation*}
\begin{aligned}
0\equiv \overline{\bar{\lambda}d\lambda+L\bar{\alpha}^{1}+\mu\alpha^{1}}
& =
\lambda d\bar{\lambda}+\bar{L}\alpha^{1}+\bar{\mu}\bar{\alpha}^{1}
=
-\bar{\lambda}d\lambda+\bar{L}\alpha^{1}+\bar{\mu}\bar{\alpha}^{1}\\
&= (L\bar{\alpha}^{1}+\mu\alpha^{1})+\bar{L}\alpha^{1}+\bar{\mu}\bar{\alpha}^{1}\\
&= (L+\bar{\mu})\bar{\alpha}^{1}+(\mu+\bar{L})\alpha^{1},
\end{aligned}
\end{equation*}
and hence by direct inspection, $\mu=-\varphi^{*}\bar{L}$. This $1$-form
\[
\tau:=\bar{\lambda} d\lambda+L\bar{\alpha}^{1}-\bar{L}\alpha^{1},
\]
which is purely imaginary, therefore vanishes upon pullback by $\varphi$. Due to the property that $\bar{\tau}=-\tau$, for any holomorphic disk $\psi:\mathbb{D}\rightarrow M^{9}$ with $\psi^{*}\omega^{0}\equiv \psi^{*}\omega^{2}\equiv \psi^{*}\bar{\omega}^{2}\equiv 0$ and $\psi^{*}(\alpha^{1}\wedge\bar{\alpha}^{1})\neq 0$, the pullback $\psi^{*}\tau$ vanishes. It is necessary to add $\tau$ to $\mathcal{I}$ and require that 
\[
d\tau\equiv 0\qquad \mod\ \mathcal{I}+\langle \tau\rangle,
\]
effectively solving for $\lambda$.

The discussion above can be related to the Chern-Moser theory, and a brief summary will be given here. Using the Appendix in Chern-Moser \cite{Chern-Moser-1974},  the most general structure  equation used is given in the following form:
\begin{equation}\label{Chern-Moser-structural-eq}
\begin{aligned}
d\alpha^{0} &= \isqrt g_{\gamma\bar{\beta}}\ \alpha^{\gamma}\wedge\alpha^{\bar{\beta}}+\alpha^{0}\wedge\phi,\\
d\alpha^{\gamma} &= \alpha^{\beta}\wedge\phi_{\beta}^{\gamma}+\alpha^{0}\wedge\phi^{\gamma},\\
d\phi &= \isqrt\alpha_{\bar{\beta}}\wedge\phi^{\bar{\beta}}+\isqrt \phi_{\bar{\beta}}\wedge\alpha^{\bar{\beta}}+\alpha^{0}\wedge\psi,\\
d\phi_{\beta}^{\gamma} &= \phi_{\beta}^{\sigma}\wedge\phi_{\sigma}^{\gamma}+\isqrt\alpha_{\beta}\wedge\phi^{\gamma}-\isqrt\phi_{\beta}\wedge\alpha^{\gamma}-\isqrt\delta_{\beta}^{\gamma}(\phi_{\sigma}\wedge\alpha^{\sigma})\\
&\hspace{0.5cm} -\frac{1}{2}\delta_{\beta}^{\gamma}\psi\wedge\alpha^{0}+\Phi_{\beta}^{\gamma},\\
d\phi^{\gamma} &= \phi\wedge\phi^{\gamma}+\phi^{\beta}\wedge\phi_{\beta}^{\gamma}-\frac{1}{2}\psi\wedge\alpha^{\gamma}+\Phi^{\gamma},\\
d\psi &= \phi\wedge\psi+2\isqrt\phi^{\beta}\wedge\phi_{\beta}+\Psi.
\end{aligned}
\end{equation}
Here the Einstein convention for summation is adopted. The 1-forms $\phi_{\bullet\bullet}$ satisfy 
\begin{equation}\label{Chern-Moser-phi}
\phi_{\gamma\bar{\beta}}+\phi_{\bar{\beta}\gamma}-g_{\gamma\bar{\beta}}\phi=0.
\end{equation}
The 2-forms that are crucial to the study of CR geometry are 
\[
\Phi_{\alpha\bar{\rho}}=S_{\alpha\beta\bar{\rho}\bar{\sigma}}\ \alpha^{\rho}\wedge\alpha^{\bar{\sigma}}+\cdots
\]
giving rise to the well-known S-tensor appearing in the expansion.

In the case of Lorentzian 5-dimensional real hypersurfaces, the matrix $g_{\alpha\bar{\beta}}$ may be chosen to be
\[
g_{\alpha\bar{\beta}}=
\left(
\begin{matrix}
1 & 0\\
0 & -1
\end{matrix}
\right).
\]
Then the first line in equation \eqref{Chern-Moser-structural-eq} reads as
\[
d\alpha^{0}=\isqrt\big(\alpha^{1}\wedge\alpha^{\bar{1}}-\alpha^{2}\wedge\alpha^{\bar{2}}\big)+\alpha^{0}\wedge\phi
\]
which is recognised as the Levi form written using the adapted basis. The $1$-form $\tau$ can alternatively be expressed in terms of $\phi_{\bullet}^{\bullet}$ as 
\[
\tau=\isqrt d\theta+e^{-\isqrt\theta}\phi_{1}^{2}+(\phi_{2}^{2}-\phi_{1}^{1})-e^{\isqrt\theta}\phi_{2}^{1}.
\]
The relations \eqref{Chern-Moser-phi} imply that $\tau$ is purely imaginary, and hence it has to be added to the ideal $\mathcal{I}$. Taking the exterior differentiation of $\tau$, the components of the $S$-tensor appears in the 2-form $d\tau$:
\[
d\tau
\equiv
-\big(
S_{22\bar{1}\bar{1}}\lambda^{2}
+4S_{21\bar{1}\bar{1}}\lambda
+6S_{11\bar{1}\bar{1}}
+4\overline{S_{21\bar{1}\bar{1}}}\bar{\lambda}
+\overline{S_{22\bar{1}\bar{1}}}\bar{\lambda}^{2}\big)\ \alpha^{1}\wedge\bar{\alpha}^{1}\qquad  \mod\ \mathcal{I}+\langle\tau\rangle.
\]
There are however differences in the case of $M^{9}\subseteq \mathbb{C}^{5}$ with Levi signature $(2,2)$. These differences come from the fact that instead of dealing with $1$-forms such as $\tau$, we are dealing with matrices of $1$-forms, as will be explained in the next subsection 1.4.

\subsection{The $(2,2)$ case and main result}

Reiterating Cartan's method, the Levi form may be written in terms of the adapted co-frame as 
\[
d\theta\equiv 
\isqrt\big(\alpha^{1}\wedge\overline{\alpha}^{1}
+
\alpha^{2}\wedge\overline{\alpha}^{2}
-
\alpha^{3}\wedge\overline{\alpha}^{3}
-
\alpha^{4}\wedge\overline{\alpha}^{4}
\big)
\qquad
\mod\ 
\theta.
\]
Assuming that a holomorphic immersion of the bi-disks $\varphi:\mathbb{D}^{2}\rightarrow M^{9}$ exists, then there is a lift $\tilde{\varphi}:\mathbb{D}^{2}\rightarrow M^{9}\times U(2)$ so that for $(s,t)\in \mathbb{D}^{2}$:
\[
\tilde{\varphi}(s,t,\bar{s},\bar{t})
=
\bigg(\varphi(s,t),
\left(
\begin{matrix}
{\sf P} & {\sf Q}\\
{\sf R} & {\sf S}
\end{matrix}
\right)
\bigg)
\eqno
{ \left(
	\begin{matrix}
	{\sf P} & {\sf Q}\\
	{\sf R} & {\sf S}
	\end{matrix}
	\right)\in U(2).}
\]
The following Pfaffian system is therefore set up:
\begin{equation}
\begin{aligned}
\omega^{0} &:= \theta,\\
\omega^{1} &:= \alpha^{1},\\
\omega^{2} &:= \alpha^{2},\\
\omega^{3} &:= \alpha^{3}-{\sf P}\alpha^{1}-{\sf Q}\alpha^{2},\\
\omega^{4} &:= \alpha^{4}-{\sf R}\alpha^{1}-{\sf S}\alpha^{2}.
\end{aligned}
\end{equation}
Since $\varphi^{*}\omega^{0}\equiv \varphi^{*}\omega^{3}\equiv\varphi^{*}\omega^{4}\equiv 0$ (along with their conjugates),  it makes sense to let $\mathcal{I}$ be the ideal generated by $\omega^{0}$, $\omega^{3}$, $\omega^{4}$, $\bar{\omega}^{3}$ and $\bar{\omega}^{4}$, describing the $\varphi_{*}\mathbb{C}T\mathbb{D}^{2}$-bundle over $\varphi(\mathbb{D}^{2})$. Their Poincar\'{e} derivatives also vanish upon pullbacks
\[
\varphi^{*}d\omega^{3}\equiv \varphi^{*}d\omega^{4}\equiv 0\qquad \mod\ \mathcal{I}.
\]
There exist complex-valued functions $A$,\dots,$J$ on $M\times U(2)$ such that
\begin{equation*}
\begin{aligned}
\left(
\begin{matrix}
{\sf \bar{P}} & {\sf \bar{R}}\\
{\sf \bar{Q}} & {\sf \bar{S}}
\end{matrix}
\right)
\left(
\begin{matrix}
d\omega^{3}\\
d\omega^{4}
\end{matrix}
\right)
&\equiv 
-
\left(
\begin{matrix}
{\sf \bar{P}} & {\sf \bar{R}}\\
{\sf \bar{Q}} & {\sf \bar{S}}
\end{matrix}
\right)
\left(
\begin{matrix}
d{\sf P} & d{\sf Q}\\
d{\sf R} & d{\sf S}
\end{matrix}
\right)
\left(
\begin{matrix}
\alpha^{1}\\
\alpha^{2}
\end{matrix}
\right)\\
& \hspace{0.5cm}
-\left(
\begin{matrix}
A\ \alpha^{1}\wedge\overline{\alpha}^{1}
+
B\ \alpha^{1}\wedge\overline{\alpha}^{2}
+
C\ \alpha^{2}\wedge\overline{\alpha}^{1}
+
D\ \alpha^{2}\wedge\overline{\alpha}^{2}
+
E\ \alpha^{1}\wedge\alpha^{2}\\
F\ \alpha^{1}\wedge\overline{\alpha}^{1}
+
G\ \alpha^{1}\wedge\overline{\alpha}^{2}
+
H\ \alpha^{2}\wedge\overline{\alpha}^{1}
+
I\ \alpha^{2}\wedge\overline{\alpha}^{2}
+
J\ \alpha^{1}\wedge\alpha^{2}
\end{matrix}
\right)\,\mod\ \mathcal{I}\\
&\equiv 0\,\mod\ \mathcal{I}.
\end{aligned}
\end{equation*}
When the torsion
\[
\left(
\begin{matrix}
A\ \alpha^{1}\wedge\overline{\alpha}^{1}
+
B\ \alpha^{1}\wedge\overline{\alpha}^{2}
+
C\ \alpha^{2}\wedge\overline{\alpha}^{1}
+
D\ \alpha^{2}\wedge\overline{\alpha}^{2}
+
E\ \alpha^{1}\wedge\alpha^{2}\\
F\ \alpha^{1}\wedge\overline{\alpha}^{1}
+
G\ \alpha^{1}\wedge\overline{\alpha}^{2}
+
H\ \alpha^{2}\wedge\overline{\alpha}^{1}
+
I\ \alpha^{2}\wedge\overline{\alpha}^{2}
+
J\ \alpha^{1}\wedge\alpha^{2}
\end{matrix}
\right)
\]
is absorbed into the first term, the equation is simplified to
\[
\left(
\begin{matrix}
{\sf \bar{P}} & {\sf \bar{R}}\\
{\sf \bar{Q}} & {\sf \bar{S}}
\end{matrix}
\right)
\left(
\begin{matrix}
d\omega^{3}\\
d\omega^{4}
\end{matrix}
\right)
=
\left(
\begin{matrix}
{\sf M_{31}}\wedge\alpha^{1} + {\sf M_{32}}\wedge\alpha^{2}\\
{\sf M_{41}}\wedge\alpha^{1} + {\sf M_{42}}\wedge\alpha^{2}
\end{matrix}
\right)\equiv 0\qquad \mod\ \mathcal{I},\]
which vanishes after being pulled back to the bi-disk. By Cartan's lemma,  a suitable modification of the $1$-forms $\varphi^{*}{\sf M_{ij}}$ by linear combinations of $\varphi^{*}\alpha^{1}$ and $\varphi^{*}\alpha^{2}$ with coefficients $l$, $m$, $n$, $p$, $q$ and $r$ is needed such that the following matrix vanishes identically on $\mathbb{D}^{2}$:
\begin{equation}\label{intro-modified-MC}
\left(
\begin{matrix}
\varphi^{*}{\sf M_{31}}-k\varphi^{*}\alpha^{1}-l\varphi^{*}\alpha^{2} & \varphi^{*} {\sf M_{32}}-l\varphi^{*}\alpha^{1}-m\varphi^{*}\alpha^{2}\\
\varphi^{*}{\sf M_{41}}-p\varphi^{*}\alpha^{1}-q\varphi^{*}\alpha^{2} & \varphi^{*}{\sf M_{42}}-q\varphi^{*}\alpha^{1}-r\varphi^{*}\alpha^{2}
\end{matrix}
\right)\equiv 0. 
\end{equation}
The main point in this paper is that the Maurer-Cartan form
\[
\left(
\begin{matrix}
{\sf \bar{P}} & {\sf \bar{R}}\\
{\sf \bar{Q}} & {\sf \bar{S}}
\end{matrix}
\right)
\left(
\begin{matrix}
d{\sf P} & d{\sf Q}\\
d{\sf R} & d{\sf S}
\end{matrix}
\right)
\]
is skew-hermitian, and so the matrix in equation \eqref{intro-modified-MC} must also be skew-hermitian, leading to a set of equations involving  $k$, $l$, $m$, $p$, $q$, $r$ that need to be solved. These equations have solutions if the following functions
\[
{\sf T_{1}} := \bar{B}+E-\bar{F},
\qquad
{\sf T_{2}} :=\bar{D}+J-\bar{H}
\]
vanish on $\tilde{\varphi}(\mathbb{D}^{2})$. The following theorem summarises the new result in this paper:

\begin{Theorem}
	Let $M^{9}\subset\mathbb{C}^{5}$ be a CR generic smooth real hypersurface passing through the origin, and whose Levi form has signature of $(2,2)$ at each point in $M^{9}$. Suppose $\varphi:\mathbb{D}^{2}\rightarrow M^{9}$ is a holomorphic immersion of bi-disk into $M^{9}$ with  $\varphi(0)=0$, then its unique lift $\tilde{\varphi}:\mathbb{D}^{2}\rightarrow M^{9}\times U(2)$ has the image lying in the zero set of two complex valued functions ${\sf T_{1}}$ and ${\sf T_{2}}$ on $M^{9}\times U(2)$.
\end{Theorem}

The existence of two complex-valued obstructions (or 4 real-valued obstructions) is the main difference from the Lorentzian case. This is to be expected because in the Lorentzian case, the Maurer-Cartan form $\bar{\lambda}d\lambda$ is purely imaginary. During the process of absorption and prolongation, only a simple modification is required for the modified Maurer-Cartan form to remain purely imaginary. Analogously, in the case where the signature is $(2,2)$, the Maurer-Cartan form is a skew-hermitian matrix. However, for most of the time, the resulting matrix after absorption of the torsions might fail to be skew-hermitian. This imposes restrictions on the torsions, resulting in the existence of 2 complex valued functions on $M\times U(2)$ which are unavoidable.

Assuming that there are no such obstructions, the Cartan process is carried on as in the case of Lorentzian CR manifolds by finding an integral element. It comes in the form  of a $2\times 2$ skew-hermitian matrix of 1-forms $\tau$  which is given by 
\[
\tau=
\left(
\begin{matrix}
-\bar{\sf P}d{\sf P}
-\bar{\sf R} d{\sf R}
&
-\bar{\sf P}d{\sf Q}
-\bar{\sf R}d{\sf S}\\
-\bar{\sf Q}d{\sf P}
-\bar{\sf S}d{\sf R}
&
-\bar{\sf Q}d{\sf Q}
-\bar{\sf S}d{\sf S}
\end{matrix}
\right)
+
\underbrace{
\left(
\begin{matrix}
\bar{A}\alpha^{1}
+\bar{B}\alpha^{2}
-A\bar{\alpha}^{1}
-B\bar{\alpha}^{2} &
\bar{F}\alpha^{1}
+\bar{G}\alpha^{2}
-C\bar{\alpha}^{1}
-D\bar{\alpha}^{2}\\
\bar{C}\alpha^{1}
+\bar{D}\alpha^{2}
-F\bar{\alpha}^{1}
-G\bar{\alpha}^{2}
&
\bar{H}\alpha^{1}
+\bar{I}\alpha^{2}
-H\bar{\alpha}^{1}
-I\bar{\alpha}^{2}
\end{matrix}
\right)}_{:=\Sigma}
\]
satisfying the following structural equation:
\[
d\tau =d\Sigma +\Sigma\wedge\Sigma\equiv 0 \qquad \mod\ \mathcal{I}+\langle \tau\rangle.
\]
In view of  Chern-Moser theory, $d\tau$ may also be expressed in terms of components of the the $S$-tensors $S_{\alpha\rho}^{\beta\sigma}$.

We conclude with the discussion of the two obstructions to holomorphic immersions of the bi-disk. By studying an example, we know that there is a certain CR smooth real-analytic hypersurface $M^{9}$ that does not satisfy the vanishing of two functions on $M^{9}\times U(2)$, concluding that such an immersion cannot happen.

\medskip\noindent 
\textbf{Acknowledgement:} The first author would like to thank Professor Robert Bryant for his explanation of Chern-Moser theory applied to the study of holomorphic disks in Lorentzian real hypersurfaces. The first author is funded by Hua Loo-Keng Center for Mathematical Sciences, AMSS, CAS, under NSFC grant number 11688101.

Grateful thanks are addressed to two anonymous established `connoisseurs' of Cartan's universe for insightful comments.


\Section{\bf The Geometric Setup of Coframes}
\label{geometric-setup}

\subsection{Notations}
The following notations will be used in this paper:
\begin{equation*}
\begin{aligned}
M^{9}:& \text{ a real hypersurface in } \mathbb{C}^{5},\\
\theta:& \text{ the contact 1-form on }M^{9},\\
\mathbb{D}^{2}:& \text{ a bidisk }\mathbb{D}\times\mathbb{D},\\
(s,t):& \text{ holomorphic variables in }\mathbb{D}^{2},\\
(z_{1},\dots,z_{5}):& \text{ holomorphic variables in } \mathbb{C}^{5}.
\end{aligned}
\end{equation*}

\subsubsection{}  Let $M^{9}$ be a real-hypersurface in $\mathbb{C}^{5}$ passing through the origin, and whose signature of the Levi form is $(2,2)$ at each point in $M^{9}$. This paper will be devoted to the study of the presence of a $2-$dimensional complex variety contained in $M^{9}$ and passing through the origin.
\begin{Definition}\label{def-holo-immersion}
A {\sl holomorphic immersion} of the bidisk $\mathbb{D}^{2}$ into $M^{9}$ is a map
\[
\varphi:\mathbb{D}^{2}\rightarrow M^{9},
\]
which is 1) an immersion of manifolds, with the property that at each point $y\in \mathbb{D}^{2}$, 2) the vector space $\varphi_{*}T_{y}\mathbb{D}^{2}$ is invariant under the complex structure $J$ and 3) there is an inclusion of the vector bundles
\[T\mathbb{D}^{2}\subseteq T^{c}M=\ker\theta.\]
\end{Definition}

\subsubsection{} The condition of $\varphi_{*}T_{y}\mathbb{D}^{2}$ being invariant under the complex structure $J$ is crucial. Without this condition, it is possible that there is a map $\varphi:\mathbb{D}^{2}\rightarrow M^{9}$ satisfying only 1) and 3) of the definition, but $\varphi$ is not holomorphic. 

\subsubsection{} To illustrate this point, let $(z_{1},z_{2},w)$ be holomorphic coordinates in $\mathbb{C}^{3}$ and consider the following hypersurface $M^{5}\subset \mathbb{C}^{3}$ given by 
\[
M^{5}:=\{w+\bar{w}=0\}.
\]
If $\varphi:\mathbb{D}\rightarrow M^{5}$ is a curve defined by $t\mapsto \varphi(t):=(t,\bar{t},0)$, then it is clearly not holomorphic. The common tangent bundle $T^{c}M^{5}$ of this real hypersurface is given at each point $x\in M^{5}$ by the following real vector space
\[
T_{x}^{c}M^{5}=
\text{span}_{\mathbb{R}}\{\partial_{z_{1}},\partial_{z_{2}},\partial_{\bar{z}_{1}},\partial_{\bar{z}_{2}}\}.
\]
On the other hand, the pushforward of $T\mathbb{D}$ by $\varphi$ is given at each point $y\in\mathbb{D}$ by 
\[
\varphi_{*}(T_{y}\mathbb{D})
=
\text{span}_{\mathbb{R}}\{\varphi_{*}\partial_{t},\varphi_{*}\partial_{\bar{t}}\}
=
\text{span}_{\mathbb{R}}\{\partial_{z_{1}}+\partial_{\bar{z}_{2}},\ \partial_{\bar{z}_{1}}+\partial_{z_{2}}\}.
\]
Therefore, it is clear that $\varphi$ is an immersion, and is contained in $T_{\varphi(y)}^{c}M^{5}$. But if $J$ is the complex structure, then $J\partial_{z_{1}}+J\partial_{\bar{z}_{2}}$ lies out of $\varphi_{*}T_{y}\mathbb{D}^{2}$ and hence $\varphi_{*}T_{y}\mathbb{D}$ is not invariant under $J$. Hence leaving out the second condition leads to an immersion which is not necessarily holomorphic.

\subsubsection{} Let $\{\alpha^{1},\alpha^{2},\alpha^{3},\alpha^{4}\}$ be a $T^{1,0*}M$ coframe which diagonalises the Levi form so that 
\begin{equation}\label{eq-geomsetup-levi-diag}
d\theta
=
\isqrt (\alpha^{1}\wedge\overline{\alpha}^{1}
+
\alpha^{2}\wedge\overline{\alpha}^{2}
-
\alpha^{3}\wedge\overline{\alpha}^{3}
-
\alpha^{4}\wedge\overline{\alpha}^{4})
\qquad
\mod\ \theta.
\end{equation}
Its respective duals $\{\mathcal{A}_{1}, \mathcal{A}_{2}, \mathcal{A}_{3}, \mathcal{A}_{4}\}$ is a $T^{1,0}M$ frame consisting of vector fields on $M^{9}$ with
\[
\alpha^{i}(\mathcal{A}_{j})=\delta_{ij},\qquad
\bar{\alpha}^{i}(\mathcal{A}_{j})=\alpha^{i}(\bar{\mathcal{A}}_{j})=0.
\]

\subsubsection{Tangent vectors of $\varphi(\mathbb{D}\times \mathbb{D})$} At each $y\in\mathbb{D}^{2}$, every vector  $v\in T_{y}^{1,0}\mathbb{D}^{2}\oplus T_{y}^{0,1}\mathbb{D}^{2}$ has its pushforward image $\varphi_{*}v$ lying in  $T_{\varphi(y)}^{1,0}M^{9}\oplus T_{\varphi(y)}^{0,1}M^{9}$ by definition \ref{def-holo-immersion}, and hence it lies in the kernel of the $1$-form
\[
\theta_{\varphi(y)}(\varphi_{*}v)=0
\eqno
({\scriptstyle v\in T_{y}^{1,0}\mathbb{D}^{2}\oplus T_{y}^{0,1}\mathbb{D}^{2}}),
\]
from which it may be inferred that the pull-back of the contact form $\theta$ by $\varphi$ vanishes identically on $\mathbb{D}^{2}$:
\[
\varphi^{*}\theta \equiv 0.
\]
Applying the Poincar\'{e} exterior differentiation of $\varphi^{*}\theta$ to both sides of the equation above, while observing that the differential operator $d$ commutes with pull-back, the Levi form $d\theta$ vanishes identically along $\varphi(\mathbb{D}^{2})$:
\[
0\equiv d\big(\varphi^{*}\theta\big)=\varphi^{*}d\theta.
\]
Taking $\{\partial_{s},\partial_{t}\}$ to be the canonical frame of the $T^{1,0}\mathbb{D}^{2}$ bundle, the identical vanishing of $d\theta$ along $\varphi(\mathbb{D}^{2})$  conducts trivially to
\[
0 = \varphi^{*}d\theta\big(\partial_{s}\wedge\partial_{\bar{s}}\big)
=
\varphi^{*}d\theta\big(\partial_{s}\wedge\partial_{\bar{t}}\big)
=
\varphi^{*}d\theta\big(\partial_{t}\wedge\partial_{\bar{s}}\big)
=
\varphi^{*}d\theta\big(\partial_{t}\wedge\partial_{\bar{t}}\big).
\]
If $\mathcal{L}_{1}:=\varphi_{*}\partial_{s}$ and $\mathcal{L}_{2}:=\varphi_{*}\partial_{t}$ are two vector fields that are tangent along $\varphi(\mathbb{D}^{2})$, the definition of the pullback of the differential forms therefore expresses the above equations in terms of $\mathcal{L}_{1}$ and $\mathcal{L}_{2}$:
\begin{equation} \label{eq-levi-vanishing-2}
\begin{aligned}
d\theta_{\varphi(y)}(\mathcal{L}_{1}|_{\varphi(y)}\wedge\overline{\mathcal{L}}_{1}|_{\varphi(y)})
=
\varphi^{*}d\theta_{y}(\partial_{s}|_{y}\wedge\partial_{\bar{s}}|_{y}) &= 0,\\
d\theta_{\varphi(y)}(\mathcal{L}_{1}|_{\varphi(y)}\wedge\overline{\mathcal{L}}_{2}|_{\varphi(y)})
=
\varphi^{*}d\theta_{y}(\partial_{s}|_{y}\wedge\partial_{\bar{t}}|_{y}) &= 0,\\
d\theta_{\varphi(y)}(\mathcal{L}_{2}|_{\varphi(y)}\wedge\overline{\mathcal{L}}_{1}|_{\varphi(y)})
=
\varphi^{*}d\theta_{y}(\partial_{t}|_{y}\wedge\partial_{\bar{s}}|_{y}) &= 0,\\
d\theta_{\varphi(y)}(\mathcal{L}_{2}|_{\varphi(y)}\wedge\overline{\mathcal{L}}_{2}|_{\varphi(y)})
=
\varphi^{*}d\theta_{y}(\partial_{t}|_{y}\wedge\partial_{\bar{t}}|_{y}) &= 0.
\end{aligned}
\end{equation}
These tangent vector fields $\mathcal{L}_{1}$ and $\mathcal{L}_{2}$ can be written as linear combinations of the frames $\mathcal{A}_{i}$ so that for some functions $f_{i}(s,t,\bar{s},\bar{t})$ and $g_{i}(s,t,\bar{s},\bar{t})$ in $\varphi(\mathbb{D}^{2})$ parametrised by $(s,t)$:
\begin{equation} \label{eq-levi-vanishing-3}
\begin{aligned}
\mathcal{L}_{1} &= \sum_{i=1}^{4}f_{i}(s,t,\bar{s},\bar{t})\mathcal{A}_{i},\\
\mathcal{L}_{2} &= \sum_{i=1}^{4}g_{i}(s,t,\bar{s},\bar{t})\mathcal{A}_{i}.
\end{aligned}
\end{equation}
Substitute each $\mathcal{L}_{i}$ in  \eqref{eq-levi-vanishing-2} by the expressions in equation \eqref{eq-levi-vanishing-3}. Using equation \eqref{eq-geomsetup-levi-diag}, the following set of relations between $f_{i}$ and $g_{i}$ is therefore obtained:
\begin{equation*}
\begin{aligned}
0 &= d\theta(\mathcal{L}_{1}\wedge\overline{\mathcal{L}}_{1}) = |f_{1}|^{2}+|f_{2}|^{2}-|f_{3}|^{2}-|f_{4}|^{2},\\
0 &= d\theta(\mathcal{L}_{1}\wedge\overline{\mathcal{L}}_{2}) = f_{1}\bar{g}_{1}+f_{2}\bar{g}_{2}-f_{3}\bar{g}_{3}-f_{4}\bar{g}_{4},\\
0 &= d\theta(\mathcal{L}_{2}\wedge\overline{\mathcal{L}}_{1}) = g_{1}\bar{f}_{1}+g_{2}\bar{f}_{2}-g_{2}\bar{f}_{3}-g_{4}\bar{f}_{4},\\
0 &= d\theta(\mathcal{L}_{2}\wedge\overline{\mathcal{L}}_{2}) = |g_{1}|^{2}+|g_{2}|^{2}-|g_{3}|^{2}-|g_{4}|^{2}.
\end{aligned}
\end{equation*}
This in turn implies that at each point $y$ in $M$ the quadratic form given by the Levi form
\[
v\mapsto d\theta_{y}(v\wedge \bar{v})
\eqno
\big({\scriptstyle 
v\in T_{y}^{1,0}M}\big)
\]
 vanishes on the $2$-dimensional plane generated by $\mathcal{L}_{1}|_{\varphi(y)}$ and $\mathcal{L}_{2}|_{\varphi(y)}$. More precisely, for any vector $(\lambda_{1},\lambda_{2})$ in $\mathbb{C}^{2}$, 
\[
d\theta
\bigg(\big(\lambda_{1}\mathcal{L}_{1}|_{\varphi(y)}+\lambda_{2}\mathcal{L}_{2}|_{\varphi(y)}\big)\wedge \overline{\big(\lambda_{1}\mathcal{L}_{1}|_{\varphi(y)}+\lambda_{2}\mathcal{L}_{2}|_{\varphi(y)}\big)}\bigg)=0.
\]
from which it can be deduced that 
\[
|\lambda_{1}f_{1}+\lambda_{2}g_{1}|^{2}
+
|\lambda_{1}f_{2}+\lambda_{2}g_{2}|^{2}
-
|\lambda_{1}f_{3}+\lambda_{2}g_{3}|^{2}
-
|\lambda_{1}f_{4}+\lambda_{2}g_{4}|^{2}
=
0.
\]
The vector subspace of $T_{y}^{1,0}M$ spanned by $\mathcal{L}_{1}|_{y}$ and $\mathcal{L}_{2}|_{y}$ is therefore in the zero set of the quadratic form
\[
(x_{1},x_{2},x_{3},x_{4})\mapsto 
x_{1}\overline{x_{1}}+x_{2}\overline{x_{2}}-
x_{3}\overline{x_{3}}-x_{4}\overline{x_{4}}.
\]

\begin{Theorem}[c.f. Sommer \cite{Sommer-1959}, proved in manuscript \cite{Merker-Manuscript-Lorentzian}] \label{Sommer's-theorem} For $n\in \mathbb{N}_{\geqslant 1}$, let $Q:\mathbb{C}^{n}\rightarrow \mathbb{R}$ be the quadratic form
\[
Q:\ (x_{1},\dots,x_{n_{+}},y_{1},\dots,y_{n_{-}})
\longmapsto 
x_{1}\overline{x_{1}}
+
\cdots
+
x_{n_{+}}\overline{x_{n_{+}}}
-
y_{1}\overline{y_{1}}
-
\cdots
-
y_{n_{-}}\overline{y_{n_{-}}}.
\]
Assume that $n_{+}\leqslant n_{-}$, and let $S$ be the set 
\[
S:=\{H\in Gr(n_+,\mathbb{C}^{n}):\ Q|_{H}\equiv 0\}.
\]
For each $H\in S$, let $\{v_{1},\dots,v_{n_+}\}$ be one of its ordered basis which can be written as
\begin{equation}
\begin{aligned}
v_{1} &= (v_{1,1},\dots,v_{n_+,1},w_{1,1},\dots,w_{n_-,1})\\
 & \hspace{0,2cm}\vdots\\
 v_{n_+} &= (v_{1,n_+},\dots,v_{n_+,n_+},w_{1,n_+},\dots,w_{n_-,n_+}).
\end{aligned}
\end{equation}
Set
\begin{equation}
B_{H}:=
\left(
\begin{matrix}
v_{1,1} & \cdots & v_{1,n_+}\\
\vdots & \ddots & \vdots\\
v_{n_+,1} & \cdots & v_{n_+,n_+}
\end{matrix}
\right)
\qquad
\text{and}
\qquad
A_{H}
:=
\left(
\begin{matrix}
w_{1,1} & \cdots & w_{1,n_+}\\
\vdots & \ddots & \vdots \\
w_{n_-,1} & \cdots & w_{n_-,n_+}
\end{matrix}
\right).
\end{equation}
Then $B_{H}$ is invertible and $A_{H}B_{H}^{-1}\in U_{n_-,n_+}$. The map
\begin{equation}
\begin{aligned}
\psi:\ S & \longrightarrow U_{n_{-},n_{+}}\\
H &\longmapsto A_{H}B_{H}^{-1}
\end{aligned}
\end{equation}
is well-defined, independent of the choice of basis of $H$. Moreover, $\psi$ describes a $1$-to-$1$ correspondence between $S$ and $U_{n_{-},n_{+}}$.
\end{Theorem}

At each point $y\in\mathbb{D}$, the plane $\text{span}_{\mathbb{C}}\{\mathcal{L}_{1}|_{\varphi(y)},\mathcal{L}_{2}|_{\varphi(y)}\}$ lies in the zero set of the quadratic equation $Q$. By the theorem, the matrix 
\[
\left(
\begin{matrix}
f_{1} & g_{1}\\
f_{2} & g_{2}
\end{matrix}
\right)
\]
is invertible, and 
there exists a unique section of the $U(2)$ bundle over $\varphi(\mathbb{D}^{2})$ which can be explicitly expressed as
\[
\left(
\begin{matrix}
\sf{p} & \sf{q}\\
\sf{r} & \sf{s}
\end{matrix}
\right)
:=
\left(
\begin{matrix}
f_{3} & g_{3}\\
f_{4} & g_{4}
\end{matrix}
\right)
\left(
\begin{matrix}
f_{1} & g_{1}\\
f_{2} & g_{2}
\end{matrix}
\right)^{-1}
\in
\Gamma(\varphi(\mathbb{D}^{2}),U(2))
\]
such that very clearly,
\begin{equation}\label{geom-setup-eq-U(2)-relation}
\left(
\begin{matrix}
f_{3}\\
f_{4}
\end{matrix}
\right)
=
\left(
\begin{matrix}
\sf{p} & \sf{q}\\
\sf{r} & \sf{s}
\end{matrix}
\right)
\left(
\begin{matrix}
f_{1}\\
f_{2}
\end{matrix}
\right)
\qquad
\text{and}
\qquad
\left(
\begin{matrix}
g_{3}\\
g_{4}
\end{matrix}
\right)
=
\left(
\begin{matrix}
\sf{p} & \sf{q}\\
\sf{r} & \sf{s}
\end{matrix}
\right)
\left(
\begin{matrix}
g_{1}\\
g_{2}
\end{matrix}
\right).
\end{equation}
Moreover, the coefficients $\sf{p}$, $\sf{q}$, $\sf{r}$, $\sf{s}$ are real analytic whenever $\mathcal{L}_{1}$ and $\mathcal{L}_{2}$ are.

\subsection{Geometric Setup of the Pfaffian system}\label{sect:geom-setup:geom-setup}

Consider the system of $1$-forms over $\varphi(\mathbb{D}^{2})$
\begin{equation}
\begin{aligned}
\omega^{0} &:= \theta,\\
\omega^{1} &:= \alpha^{1},\\
\omega^{2} &:= \alpha^{2},\\
\omega^{3} &:= \alpha^{3}-\sf{p}\alpha^{1}-\sf{q}\alpha^{2},\\
\omega^{4} &:= \alpha^{4} - \sf{r}\alpha^{1}-\sf{s}\alpha^{2}.
\end{aligned}
\end{equation}
It is clear that the set of vectors $\{\theta,\omega^{1},\omega^{2},\omega^{3},\omega^{4}\}$ is  linearly independent at each point of $M$, and furthermore for $k=3,\ 4$,
\begin{equation*}
\begin{aligned}
\varphi^{*}\omega^{k}(\partial_{s}) &= \omega^{k}(\varphi_{*}\partial_{s})=\omega^{k}(\mathcal{L}_{1})\equiv 0,\\
\varphi^{*}\omega^{k}(\partial_{t}) &= \omega^{k}(\varphi_{*}\partial_{t})=\omega^{k}(\mathcal{L}_{2})\equiv 0,\\
\varphi^{*}\overline{\omega}^{k}(\partial_{s}) &= \overline{\omega}^{k}(\varphi_{*}\partial_{s})=\overline{\omega}^{k}(\mathcal{L}_{1})\equiv 0,\\
\varphi^{*}\overline{\omega}^{k}(\partial_{t}) &= \overline{\omega}^{k}(\varphi_{*}\partial_{t})=\overline{\omega}^{k}(\mathcal{L}_{2})\equiv 0,
\end{aligned}
\hspace{1cm}
\begin{aligned}
\varphi^{*}\omega^{k}(\partial_{\bar{s}}) &= \omega^{k}(\varphi_{*}\partial_{\bar{s}})=\omega^{k}(\overline{\mathcal{L}}_{1})\equiv 0,\\
\varphi^{*}\omega^{k}(\partial_{\bar{t}}) &= \omega^{k}(\varphi_{*}\partial_{\bar{t}})=\omega^{k}(\overline{\mathcal{L}}_{2})\equiv 0,\\
\varphi^{*}\overline{\omega}^{k}(\partial_{\bar{s}}) &= \overline{\omega}^{k}(\varphi_{*}\partial_{\bar{s}})=\overline{\omega}^{k}(\overline{\mathcal{L}}_{1})\equiv 0,\\
\varphi^{*}\overline{\omega}^{k}(\partial_{\bar{t}}) &= \overline{\omega}^{k}(\varphi_{*}\partial_{\bar{t}})=\overline{\omega}^{k}(\overline{\mathcal{L}}_{2})\equiv 0,
\end{aligned}
\end{equation*}
where the vanishing follows from equation \eqref{geom-setup-eq-U(2)-relation}. Therefore, the following $1$ forms vanish on $\mathbb{C}T\mathbb{D}^{2}$ after being pulled back to $\mathbb{D}^{2}$:
\[
\varphi^{*}\theta\equiv \varphi^{*}\omega^{3}\equiv \varphi^{*}\omega^{4}\equiv \varphi^{*}\overline{\omega}^{3}\equiv \varphi^{*}\overline{\omega}^{4}\equiv 0.
\] 
Finally, 
\begin{equation}\label{geom-setup-vol-form}
\alpha^{1}\wedge\overline{\alpha}^{1}\wedge\alpha^{2}\wedge\overline{\alpha^{2}}
\bigg(
\mathcal{L}_{1}\wedge\overline{\mathcal{L}_{1}}\wedge
\mathcal{L}_{2}\wedge\overline{\mathcal{L}_{2}}
\bigg)
=
|f_{1}g_{2}-f_{2}g_{1}|^{2},
\end{equation}
and the member of the right hand side of the equation is non-zero everywhere on $\mathbb{D}^{2}$ due to Sommer's theorem. Therefore, the $(2,2)$-form $\alpha^{1}\wedge\bar{\alpha}^{1}\wedge\alpha^{2}\wedge\bar{\alpha}^{2}$ is never zero at each point of $\varphi(\mathbb{D}^{2})$.


\Section{\bf The General Setup}
\label{general-setup}


\subsubsection{The Pfaffian System} Based on the discussion in the previous section, the following Pfaffian system is a necessary condition for a holomorphic bidisk to be immersed in $M^{9}$:
\begin{equation} \label{eq:general-1} 
\begin{aligned}
\omega^{0} &= \theta,\\
\omega^{1} &= \alpha^{1},\\
\omega^{2} &= \alpha^{2},\\
\omega^{3} &= \alpha^{3} - {\sf P}\alpha^{1}-{\sf Q} \alpha^{2},\\
\omega^{4} &= \alpha^{4} - {\sf R}\alpha^{1} - {\sf S}\alpha^{2},
\end{aligned}
\end{equation}
with the following matrix
\[
\left(
\begin{matrix}
{\sf P} & {\sf Q}\\
{\sf R} & {\sf S}
\end{matrix}
\right)
\in
U(2).
\]
Let $\mathcal{I}$ be the differential ideal
\[
\mathcal{I} := \langle \omega^{0},\ \omega^{3},\ \omega^{4}, \bar{\omega}^{3},\ \bar{\omega}^{4}\rangle,
\]
and let $ \delta:=\varphi^{*}\alpha^{1}$, $\varepsilon:=\varphi^{*}\alpha^{2}$. Based on the discussion in section \ref{sect:geom-setup:geom-setup}, the integral element
\[
(\delta\wedge\overline{\delta}\wedge\varepsilon\wedge\overline{\varepsilon})|_{y}\neq 0
\eqno
(
{\scriptstyle \forall y\in\mathbb{D}^{2}}
)
\]
is non-vanishing everywhere on $\mathbb{D}^{2}$.

\subsubsection{}
To understand the Pfaffian system through the point of view of pullbacks of differential forms, for each $1\leqslant i\leqslant 4$, let ${\sf a_{i}}$ and ${\sf b_{i}}$ be functions on $\mathbb{D}^{2}$ such that 
\[
\varphi^{*}\alpha^{i}={\sf a_{i}}(s,t,\bar{s},\bar{t})\ ds + {\sf b_{i}}(s,t,\bar{s},\bar{t})\ dt.
\]
The $(0,1)$ part of $\varphi^{*}\alpha^{i}$ vanishes identically because $\varphi$ is a holomorphic map. The vanishing of the Levi form after pullback results in 
\begin{equation}\label{eq:gensetup:pullback-1-forms}
\begin{aligned}
0 &= \varphi^{*}\alpha^{1}\wedge\varphi^{*}\overline{\alpha}^{1}
+
\varphi^{*}\alpha^{2}\wedge\varphi^{*}\overline{\alpha}^{2}
-
\varphi^{*}\alpha^{3}\wedge\varphi^{*}\overline{\alpha}^{3}
-
\varphi^{*}\alpha^{4}\wedge\varphi^{*}\overline{\alpha}^{4}\\
&=
(|{\sf a_{1}}|^{2}+|{\sf a_{2}}|^{2}-|{\sf a_{3}}|^{2}-|{\sf a_{4}}|^{2})\ ds\wedge d\bar{s}
+
({\sf a_{1}}\overline{\sf b}_{1}+{\sf a_{2}}\overline{\sf b}_{2}-{\sf a_{3}}\overline{\sf b}_{3}-{\sf a_{4}}\overline{\sf b}_{4})\ ds\wedge d\bar{t}\\
&+
({\sf b_{1}}\overline{\sf a}_{1}+{\sf b_{2}}\overline{\sf a}_{2}-{\sf b_{3}}\overline{\sf a}_{3}-{\sf b_{4}}\overline{\sf a}_{4})\ dt\wedge d\bar{s}
+
(|{\sf b_{1}}|^{2}+|{\sf b_{2}}|^{2}-|{\sf b_{3}}|^{2}-|{\sf b_{4}}|^{2})\ dt\wedge d\bar{t}.
\end{aligned}
\end{equation}
In addition, by the hypothesis that
\[
\varphi^{*}(\alpha^{1}\wedge\overline{\alpha}^{1}\wedge\alpha^{2}\wedge\overline{\alpha}^{2})
=
|{\sf a_{1}}{\sf b_{2}}-{\sf b_{1}}{\sf a_{2}}|^{2}\ ds\wedge d\bar{s}\wedge dt\wedge d\bar{t}
\neq 0
\]
at each point in $\mathbb{D}^{2}$, the tuples $({\sf a_{1}},{\sf a_{2}})$ and $({\sf b_{1}},{\sf b_{2}})$ are not colinear everywhere, and therefore the two vectors $({\sf a_{1}},{\sf a_{2}},{\sf a_{3}},{\sf a_{4}})$ and $({\sf b_{1}},{\sf b_{2}},{\sf b_{3}},{\sf b_{4}})$ are also not colinear. 

Consequently from equation \eqref{eq:gensetup:pullback-1-forms}, the distribution of rank 2 complex planes over $\varphi(\mathbb{D}^{2})$ generated by the linearly independent vectors $({\sf a_{1}},{\sf a_{2}},{\sf a_{3}},{\sf a_{4}})$ and $({\sf b_{1}},{\sf b_{2}},{\sf b_{3}},{\sf b_{4}})$ lie in the solution space of the quadratic equation $x_{1}\bar{x}_{1}+x_{2}\bar{x}_{2}-y_{1}\bar{y}_{1}-y_{2}\bar{y}_{2}=0$. Again by Sommer, there exists a unique section of the unitary bundle over $\mathbb{D}^{2}$ explicitly expressed as
\[
\left(
\begin{matrix}
{\sf p'} & {\sf q'}\\
{\sf r'} & {\sf s'}
\end{matrix}
\right)
:=
\left(
\begin{matrix}
{\sf a_{3}} & {\sf b_{3}}\\
{\sf a_{4}} & {\sf b_{4}}
\end{matrix}
\right)
\left(
\begin{matrix}
{\sf a_{1}} & {\sf b_{1}}\\
{\sf a_{2}} & {\sf b_{2}}
\end{matrix}
\right)^{-1}\in \Gamma(\mathbb{D}^{2},U(2))
\]
so that
\begin{equation}\label{eq:gensetup:p'q'r's'}
\begin{aligned}
\left(
\begin{matrix}
{\sf a_{3}}\\
{\sf a_{4}}
\end{matrix}
\right)
&=
\left(
\begin{matrix}
{\sf p'} & {\sf q'}\\
{\sf r'} & {\sf s'}
\end{matrix}
\right)
\left(
\begin{matrix}
{\sf a_{1}}\\
{\sf a_{2}}
\end{matrix}
\right),\\
\left(
\begin{matrix}
{\sf b_{3}}\\
{\sf b_{4}}
\end{matrix}
\right)
&=
\left(
\begin{matrix}
{\sf p'} & {\sf q'}\\
{\sf r'} & {\sf s'}
\end{matrix}
\right)
\left(
\begin{matrix}
{\sf b_{1}}\\
{\sf b_{2}}
\end{matrix}
\right).
\end{aligned}
\end{equation}
\begin{Proposition}
At each $y\in\mathbb{D}^{2}$,
\[
\left(
\begin{matrix}
{\sf p'} & {\sf q'}\\
{\sf r'} & {\sf s'}
\end{matrix}
\right)
=
\left(
\begin{matrix}
{\sf p} & {\sf q}\\
{\sf r} & {\sf s}
\end{matrix}
\right).
\]
\end{Proposition}
\begin{proof}
The following vanishing identities 
\begin{equation}
\begin{aligned}
0 &\equiv 
\varphi^{*}\omega^{3}
=
({\sf a_{3}}ds+{\sf b_{3}}dt)-
{\sf p}({\sf a_{1}}ds+{\sf b_{1}}dt)
-
{\sf q}({\sf a_{2}}ds+{\sf b_{2}}dt),\\
0 &\equiv 
\varphi^{*}\omega^{4}
=
({\sf a_{4}}ds+{\sf b_{4}}dt)-
{\sf r}({\sf a_{1}}ds+{\sf b_{1}}dt)
-
{\sf s}({\sf a_{2}}ds+{\sf b_{2}}dt),
\end{aligned}
\end{equation}
result in having
\begin{equation}
\begin{aligned}
\left(
\begin{matrix}
{\sf a_{3}}\\
{\sf a_{4}}
\end{matrix}
\right)
=
\left(
\begin{matrix}
{\sf p} & {\sf q}\\
{\sf r} & {\sf s}
\end{matrix}
\right)
\left(
\begin{matrix}
{\sf a_{1}}\\
{\sf a_{2}}
\end{matrix}
\right),
\qquad
\left(
\begin{matrix}
{\sf b_{3}}\\
{\sf b_{4}}
\end{matrix}
\right)
=
\left(
\begin{matrix}
{\sf p} & {\sf q}\\
{\sf r} & {\sf s}
\end{matrix}
\right)
\left(
\begin{matrix}
{\sf b_{1}}\\
{\sf b_{2}}
\end{matrix}
\right),
\end{aligned}
\end{equation}
or in the matrix form,
\[
\left(
\begin{matrix}
{\sf a_{3}} & {\sf b_{3}}\\
{\sf a_{4}} & {\sf b_{4}}
\end{matrix}
\right)
=
\left(
\begin{matrix}
{\sf p} & {\sf q}\\
{\sf r} & {\sf s}
\end{matrix}
\right)
\left(
\begin{matrix}
{\sf a_{1}} & {\sf b_{1}}\\
{\sf a_{2}} & {\sf b_{2}}
\end{matrix}
\right).
\]
Since $({\sf a_{1}},{\sf a_{2}})$ and $({\sf b_{1}},{\sf b_{2}})$ are not colinear due to the condition on the integral element $\varphi^{*}(\alpha^{1}\wedge\overline{\alpha}^{1}\wedge\alpha^{2}\wedge\overline{\alpha}^{2})\neq 0$, the matrix
\[
\left(
\begin{matrix}
{\sf a_{1}} & {\sf b_{1}}\\
{\sf a_{2}} & {\sf b_{2}}
\end{matrix}
\right)
\]
is invertible and hence
\[
\left(
\begin{matrix}
{\sf p'} & {\sf q'}\\
{\sf r'} & {\sf s'}
\end{matrix}
\right)
=
\left(
\begin{matrix}
{\sf a_{3}} & {\sf b_{3}}\\
{\sf a_{4}} & {\sf b_{4}}
\end{matrix}
\right)
\left(
\begin{matrix}
{\sf a_{1}} & {\sf b_{1}}\\
{\sf a_{2}} & {\sf b_{2}}
\end{matrix}
\right)^{-1}
=
\left(
\begin{matrix}
{\sf p} & {\sf q}\\
{\sf r} & {\sf s}
\end{matrix}
\right).\qedhere
\]
\end{proof}

\subsubsection{The unique lift $\tilde{\varphi}$} It is important to verify that whenever the Pfaffian system given by \eqref{eq:general-1} gives an immersion $\varphi:\mathbb{D}^{2}\rightarrow M^{9}$, its image $\varphi(\mathbb{D}^{2})$ is holomorphic. This means that the map $\varphi$ has to satisfy the three conditions in definition \ref{def-holo-immersion}. The immersion $\varphi:\mathbb{D}^{2}\rightarrow M^{9}$ has a unique lift  $\tilde{\varphi}:\mathbb{D}^{2}\rightarrow M^{9}\times U(2)$ with
\[
\tilde{\varphi}(y)=
\left(\varphi(y),
\left(
\begin{matrix}
{\sf f_{3}} & {\sf g_{3}}\\
{\sf f_{4}} & {\sf g_{4}}
\end{matrix}
\right)
\left(
\begin{matrix}
{\sf f_{1}} & {\sf g_{1}}\\
{\sf f_{2}} & {\sf g_{2}}
\end{matrix}
\right)^{-1}
\right)
\]
 so that the following diagram commutes
\[
\xymatrix{
 & M^{9}\times U(2) \ar[d]^{\pi}\\
\mathbb{D}^{2} \ar[r]_{\varphi}\ar[ur]^{\tilde{\varphi}} & M^{9}
}
\]
where $\pi$ is the usual projection map onto the first factor.

At this stage, it would be useful to introduce the sub-bundle of $T^{1,0}M$ in the following way. For each fixed fixed point $(x,{\sf P},{\sf Q}, {\sf R},{\sf S})\in M\times U(2)$, let 
\begin{equation}
\begin{aligned}
L_{(x,{\sf P},{\sf Q}, {\sf R},{\sf S})}^{1,0} &:=
\{
v\in T_{x}^{1,0}M:\  0=\alpha^{3}-{\sf P}\alpha^{1}-{\sf Q}\alpha^{2},\ 0=\alpha^{4}-{\sf R}\alpha^{1}-{\sf S}\alpha^{2}\},\\
L_{(x,{\sf P},{\sf Q}, {\sf R},{\sf S})}^{0,1} &:=
\{
v\in T_{x}^{0,1}M:\  0=\overline{\alpha}^{3}-\bar{{\sf P}}\overline{\alpha}^{1}-\bar{{\sf Q}}\overline{\alpha}^{2},\ 0=\overline{\alpha}^{4}-\bar{{\sf R}}\overline{\alpha}^{1}-\bar{{\sf S}}\overline{\alpha}^{2}\},
\end{aligned}
\end{equation}
and hence
\[
\{
0=\omega^{0}=\omega^{3}=\omega^{4}=\bar{\omega}^{3}=\bar{\omega}^{4}
\}
=
L_{(x,{\sf P},{\sf Q}, {\sf R},{\sf S})}^{1,0}\oplus 
L_{(x,{\sf P},{\sf Q}, {\sf R},{\sf S})}^{0,1}.
\]
The common bundle $L_{(x,{\sf P},{\sf Q}, {\sf R},{\sf S})}^{c}=\Re(L_{(x,{\sf P},{\sf Q}, {\sf R},{\sf S})}^{1,0}) =\Re(L_{(x,{\sf P},{\sf Q}, {\sf R},{\sf S})}^{0,1})$, is invariant under the complex structure $J$.

\begin{Proposition}
Let $\varphi:\mathbb{D}^{2}\rightarrow M^{9}$ be an immersion passing through the origin, and let $\tilde{\varphi}:\mathbb{D}^{2}\rightarrow M^{9}\times U(2)$ be its unique lift. Let $\pi:M^{9}\times U(2)\rightarrow M^{9}$ be the projection map onto the first factor. If $\varphi$ is a holomorphic immersion, then it satisfies 
\begin{equation}\label{eq:pfaffian-gives-hol-imm}
0 = \tilde{\varphi}^{*}(\omega^{0})=
\tilde{\varphi}^{*}\omega^{3}
=\tilde{\varphi}^{*}\omega^{4}=
\tilde{\varphi}^{*}\overline{\omega}^{3}
=
\tilde{\varphi}^{*}\overline{\omega}^{4}.
\end{equation}
Conversely, if $\tilde{\varphi}$ satisfies \eqref{eq:pfaffian-gives-hol-imm}, then the image $\varphi(\mathbb{D}^{2})$ is holomorphic.
\end{Proposition}

\begin{proof}
The direct implication follows immediately from the definition of $L_{(x,{\sf P},{\sf Q},{\sf R},{\sf S})}^{1,0}$, and equations \eqref{eq:gensetup:p'q'r's'}.

Before proving the converse, recall that $\varphi$ satisfies the immersion criteria 1 in definition \ref{def-holo-immersion}. It remains to prove the other two criteria. Criteria 3 of definition \ref{def-holo-immersion}, which is $\varphi_{*}T\mathbb{D}^{2}\subset \ker\theta$, is automatic since 
\[
\varphi^{*}\omega^{0}=0
\]
where the vanishing follows from hypothesis. 

It remains to prove invariance under the complex action $J$. From the following vanishing identities
\begin{equation}
\begin{aligned}
0 &= \tilde{\varphi}^{*}\omega^{3} = 
\varphi^{*}\alpha^{3}-{\sf p}\varphi^{*}\alpha^{1}
-
{\sf q}\varphi^{*}\alpha^{2},\\
0 &= \tilde{\varphi}^{*}\omega^{4} = 
\varphi^{*}\alpha^{4}-{\sf r}\varphi^{*}\alpha^{1}
-
{\sf s}\varphi^{*}\alpha^{2},\\
0 &= \tilde{\varphi}^{*}\overline{\omega}^{3} = 
\varphi^{*}\overline{\alpha}^{3}-{\sf \bar{p}}\varphi^{*}\overline{\alpha}^{1}
-
{\sf \bar{q}}\varphi^{*}\overline{\alpha}^{2},\\
0 &= \tilde{\varphi}^{*}\overline{\omega}^{4} = 
\varphi^{*}\overline{\alpha}^{3}-{\sf \bar{r}}\varphi^{*}\overline{\alpha}^{1}
-
{\sf \bar{s}}\varphi^{*}\overline{\alpha}^{2},
\end{aligned}
\end{equation}
it follows that \[\varphi_{*}(\mathbb{C}T_{y}\mathbb{D}^{2})\subseteq 
L_{\tilde{\varphi}(y)}^{1,0}
\oplus 
L_{\tilde{\varphi}(y)}^{0,1}.
\]
which turns out to be an equality since both have the same rank. Taking the real part,
\begin{equation}
\begin{aligned}
\varphi_{*}T_{y}\mathbb{D}^{2}
=
\varphi_{*}\Re\big(\mathbb{C}T_{y}\mathbb{D}^{2}\big)
=
\Re\varphi_{*}\mathbb{C}T_{y}\mathbb{D}^{2}
=
\Re \big(L_{\tilde{\varphi}(y)}^{1,0}
\oplus 
L_{\tilde{\varphi}(y)}^{0,1}\big)
=
L_{\tilde{\varphi}(y)}^{c}.
\end{aligned}
\end{equation}
By a theorem of Levi-Civita \cite[page 99]{Shabat-1992}, this implies that the image $\varphi(\mathbb{D}^{2})$ is a complex manifold.
\end{proof}

\subsubsection{} If $U$ is a section of the $U(2)$-bundle, then $U^{*}dU$ is a skew-hermitian matrix. This is because $U^{*}U=I$ and so taking the exterior derivative of this identity yields the following relation
\[
0=d\big(U^{*}U\big)=dU^{*}\cdot U+U^{*}dU,
\]
from which it may be easily deduced that
\[
\big(
U^{*}dU
\big)^*
=
dU^{*}\cdot U
=
-U^{*}dU.
\]
 
\subsubsection{} Adopting the matrix representation of the last two lines of the general Pfaffian setup in equation \eqref{eq:general-1},
\[
\left(
\begin{matrix}
\omega^{3}\\
\omega^{4}
\end{matrix}
\right)
=
-\left(
\begin{matrix}
{\sf P} & {\sf Q}\\
{\sf R}& {\sf S}
\end{matrix}
\right)
\left(
\begin{matrix}
\alpha^{1}\\
\alpha^{2}
\end{matrix}
\right)
+
\left(
\begin{matrix}
\alpha^{3}\\
\alpha^{4}
\end{matrix}
\right).
\]
Taking the exterior derivative on both sides of the equation above to obtain
\begin{equation}\label{eq:domega3-and-domega4}
\begin{aligned}
\left(
\begin{matrix}
d\omega^{3}\\
d\omega^{4}
\end{matrix}
\right)
&=
-\left(
\begin{matrix}
d{\sf P} & d{\sf Q}\\
d{\sf R} & d{\sf S}
\end{matrix}
\right)
\wedge
\left(
\begin{matrix}
\alpha^{1}\\
\alpha^{2}
\end{matrix}
\right)
-\left(
\begin{matrix}
{\sf P} & {\sf Q}\\
{\sf R} & {\sf S}
\end{matrix}
\right)
\wedge
\left(
\begin{matrix}
d\alpha^{1}\\
d\alpha^{2}
\end{matrix}
\right)
+
\left(
\begin{matrix}
d\alpha^{3}\\
d\alpha^{4}
\end{matrix}
\right).
\end{aligned}
\end{equation}
On the other hand, by making a transformation of coframes, followed by taking exterior derivative, and finally using the equation \eqref{eq:domega3-and-domega4},
\begin{equation*} 
\begin{aligned}
d\left(
\left(
\begin{matrix}
\bar{\sf P} & \bar{\sf R}\\
\bar{\sf Q} & \bar{\sf S}
\end{matrix}
\right)
\left(
\begin{matrix}
\omega^{3}\\
\omega^{4}
\end{matrix}
\right)
\right)
&=
\left(
\begin{matrix}
d\bar{\sf P} & d\bar{\sf R}\\
d\bar{\sf Q} & d\bar{\sf S}
\end{matrix}
\right)
\wedge
\left(
\begin{matrix}
\omega^{3}\\
\omega^{4}
\end{matrix}
\right)
+
\left(
\begin{matrix}
\bar{\sf P} & \bar{\sf R} \\
\bar{\sf Q} & \bar{\sf S}
\end{matrix}
\right)
\left(
\begin{matrix}
d\omega^{3}\\
d\omega^{4}
\end{matrix}
\right)\\
&\equiv 
\left(
\begin{matrix}
\bar{\sf P} & \bar{\sf R}\\
\bar{\sf Q} & \bar{\sf S}
\end{matrix}
\right)
\left(
\begin{matrix}
d\omega^{3}\\
d\omega^{4}
\end{matrix}
\right)
\qquad 
\mod\ \mathcal{I}\\
&\equiv
-\left(\begin{matrix}
\bar{\sf P} & \bar{\sf R}\\
\bar{\sf Q} & \bar{\sf S}
\end{matrix}
\right)
\left(
\begin{matrix}
d{\sf P} & d{\sf Q}\\
d{\sf R} & d{\sf S}
\end{matrix}
\right)
\wedge
\left(
\begin{matrix}
\alpha^{1}\\
\alpha^{2}
\end{matrix}
\right)
-
\left(
\begin{matrix}
d\alpha^{1}\\
d\alpha^{2}
\end{matrix}
\right)
+
\left(\begin{matrix}
\bar{\sf P} & \bar{\sf R}\\
\bar{\sf Q} & \bar{\sf S}
\end{matrix}
\right)
\left(
\begin{matrix}
d\alpha^{3}\\
d\alpha^{4}
\end{matrix}
\right)\,\,\mod\ \mathcal{I}.
\end{aligned}
\end{equation*}
The reason for this change of coframe is to allow the Maurer-Cartan form $U^{*}dU$ to appear. The skew-hermitian property of this Maurer-Cartan form will necessarily impose the same structure on the torsion. If the torsion does not follow the anti-hermitian property, it then constitutes an obstruction to the presence of holomorphic bi-disk in $M^{9}$ passing through the origin.

\subsubsection{} It remains to study the 2-forms $d\alpha^{i}$. Note that for $1\leqslant i\leqslant 4$, each $\alpha^{i}$ can be written as
\begin{equation}\label{changeofcoframetodiagonalise}
\alpha^{i}
=
\sum_{j=1}^{4}
P_{ij}(z,\bar{z},v)\ dz^{j},
\end{equation}
where the matrix $(P_{ij})$ which is invertiable at each point in $M^{9}$, describes the change of coframes. Taking the exterior derivative,
 \begin{equation*}
\begin{aligned}
d\alpha^{i}
&=
\sum_{j=1}^{4}
dP_{ij}\wedge dz^{j}\\
&=
\sum_{j=1}^{4}\bigg[
\sum_{k=1}^{4}
\left(
\frac{\partial P_{ij}}{\partial z_{k}}
dz^{k}
+
\frac{\partial P_{ij}}{\partial\bar{z}_{k}}
d\bar{z}^{k}
\right)
+
\frac{\partial P_{ij}}{\partial v}dv
\bigg]\wedge dz^{j}\\
&=
\sum_{1\leqslant j<k\leqslant 4}
-\left(
\frac{\partial P_{ij}}{\partial z_{k}}
-
\frac{\partial P_{ik}}{\partial z_{j}}
\right)\ dz^{j}\wedge dz^{k}
+
\sum_{j,k=1}^{4}
\frac{\partial P_{ij}}{\partial \bar{z}_{k}}\ 
d\bar{z}^{k}\wedge dz^{j}
+
\sum_{j=1}^{4}
\frac{\partial P_{ij}}{\partial v}\ dv\wedge dz^{j}.
\end{aligned}
\end{equation*}
Computing modulo $\theta$, which is
\[
\theta = -dv+\sum_{k=1}^{4}A_{k}dz^{k}+\sum_{k=1}^{4}\bar{A}_{k}d\bar{z}^{k},
\]
replaces $dv$ by $\sum A_{k}dz^{k}+\sum \bar{A}_{k}d\bar{z}^{k}$ in $d\alpha^{i}$, resulting in 
{\footnotesize
\begin{equation*}
\begin{aligned}
d\alpha^{i} &\equiv 
\sum_{1\leqslant j<k\leqslant 4}
\left(
\frac{\partial P_{ij}}{\partial z_{k}}
-
\frac{\partial P_{ik}}{\partial z_{j}}
\right)\ dz^{j}\wedge dz^{k}
+
\sum_{j,k=1}^{4}
\frac{\partial P_{ij}}{\partial \bar{z}_{k}}\ d\bar{z}^{k}\wedge dz^{j}
+
\sum_{j=1}^{4}
\frac{\partial P_{ij}}{\partial v}
\left(
\sum_{l=1}^{4}A_{l}dz^{l}+
\sum_{l=1}^{4}\bar{A}_{l}d\bar{z}^{l}
\right)\wedge dz^{j}\\
&=
\sum_{1\leqslant j<k\leqslant 4}
P_{ijk}\ dz^{i}\wedge dz^{j}
+
\sum_{j,k=1}^{4}
Q_{ijk}\ dz^{j}\wedge d\bar{z}^{k}\ \ \ \ \mod\ \theta,
\end{aligned}
\end{equation*}
}
which is a sum of $(1,1)$ and $(2,0)$ forms. Applying the inverse of the change of coframe in equation  \eqref{changeofcoframetodiagonalise}, it follows that $d\alpha^{i}$ is a sum of the two forms  $\alpha^{i}\wedge\alpha^{j}$ and $\alpha^{i}\wedge\bar{\alpha}^{j}$.
 
\subsubsection{} 
By the remark in the previous paragraph, there exist torsions $A$,\ \dots,\ $J$, such that 
\begin{equation}\label{eq:general-2}
\begin{aligned}
\left(
\begin{matrix}
\bar{\sf P} & \bar{\sf R}\\
\bar{\sf Q} & \bar{\sf S}
\end{matrix}
\right)
\left(
\begin{matrix}
d\omega^{3}\\
d\omega^{4}
\end{matrix}
\right)
&=
-\left(
\begin{matrix}
\bar{\sf P} & \bar{\sf R}\\
\bar{\sf Q} & \bar{\sf S}
\end{matrix}
\right)
\left(
\begin{matrix}
d{\sf P} & d{\sf Q}\\
d{\sf R} & d{\sf S}
\end{matrix}
\right)
\wedge
\left(
\begin{matrix}
\alpha^{1}\\
\alpha^{2}
\end{matrix}
\right)
+
\left[
\left(
\begin{matrix}
\bar{\sf P} & \bar{\sf R}\\
\bar{\sf Q} & \bar{\sf S}
\end{matrix}
\right)
\left(
\begin{matrix}
d\alpha^{3}\\
d\alpha^{4}
\end{matrix}
\right)
-
\left(
\begin{matrix}
d\alpha^{1}\\
d\alpha^{2}
\end{matrix}
\right)
\right]\\
&=
-\left(
\begin{matrix}
\bar{\sf P} & \bar{\sf R}\\
\bar{\sf Q} & \bar{\sf S}
\end{matrix}
\right)
\left(
\begin{matrix}
d{\sf P} & d{\sf Q}\\
d{\sf R} & d{\sf S}
\end{matrix}
\right)
\wedge
\left(
\begin{matrix}
\alpha^{1}\\
\alpha^{2}
\end{matrix}
\right)\\
&\hspace{0.5cm}+
\left(
\begin{matrix}
A\alpha^{1}\wedge\bar{\alpha}^{1}
+ B\alpha^{1}\wedge\bar{\alpha}^{2}
+ C\alpha^{2}\wedge\bar{\alpha}^{1}
+ D\alpha^{2}\wedge\bar{\alpha}^{2}
+ E\alpha^{1}\wedge\alpha^{2}\\
F\alpha^{1}\wedge\bar{\alpha}^{1}
+ G\alpha^{1}\wedge\bar{\alpha}^{2}
+ H\alpha^{2}\wedge\bar{\alpha}^{1}
+ I\alpha^{2}\wedge\bar{\alpha}^{2}
+ J\alpha^{1}\wedge\alpha^{2}
\end{matrix}
\right)
\qquad \mod\ \mathcal{I}.
\end{aligned}
\end{equation} 
The main objects of study in the rest of the paper are the torsions and their relations between each other.

\subsection{Conditions on Torsions} 
\subsubsection{}
Write $\tilde{\omega}^{3}$ and $\tilde{\omega}^{4}$ as 
\[
\left(
\begin{matrix}
\tilde{\omega}^{3}\\
\tilde{\omega}^{4}
\end{matrix}
\right)
:=
\left(
\begin{matrix}
\bar{\sf P} & \bar{\sf R}\\
\bar{\sf Q} & \bar{\sf S}
\end{matrix}
\right)
\left(
\begin{matrix}
\omega^{3}\\
\omega^{4}
\end{matrix}
\right),
\]
and also let
\[
\left(
\begin{matrix}
{\sf M_{31}} & {\sf M_{32}}\\
{\sf M_{41}} & {\sf M_{42}}
\end{matrix}
\right)
:=
\left(
\begin{matrix}
-\bar{\sf P}d{\sf P} - \bar{\sf R}d{\sf R} & 
-\bar{\sf P}d{\sf Q} - \bar{\sf R}d{\sf S}\\
-\bar{\sf Q}d{\sf P} - \bar{\sf S}d{\sf R} & 
-\bar{\sf Q}d{\sf Q} - \bar{\sf S}d{\sf S}
\end{matrix}
\right).
\]
From equation \eqref{eq:general-2}, 
\begin{equation}\label{eq:general-3}
\begin{aligned}
d\tilde{\omega}^{3} &\equiv
{\sf M_{31}}\wedge\alpha^{1}+{\sf M_{32}}\wedge\alpha^{2}\\ 
& \hspace{0.5cm} +A\alpha^{1}\wedge\bar{\alpha}^{1}
+ B\alpha^{1}\wedge\bar{\alpha}^{2}
+ C\alpha^{2}\wedge\bar{\alpha}^{1}
+ D\alpha^{2}\wedge\bar{\alpha}^{2}
+ E\alpha^{1}\wedge\alpha^{2}\qquad \mod\ \mathcal{I},\\
d\tilde{\omega}^{4} &\equiv
{\sf M_{41}}\wedge\alpha^{1}+{\sf M_{42}}\wedge\alpha^{2}\\
& \hspace{0.5cm} +F\alpha^{1}\wedge\bar{\alpha}^{1}
+ G\alpha^{1}\wedge\bar{\alpha}^{2}
+ H\alpha^{2}\wedge\bar{\alpha}^{1}
+ I\alpha^{2}\wedge\bar{\alpha}^{2}
+ J\alpha^{1}\wedge\alpha^{2}\qquad \mod\ \mathcal{I},
\end{aligned}
\end{equation}
where the following skew-hermitian property is observed
\[
\left(
\begin{matrix}
{\sf M_{31}} & {\sf M_{32}}\\
{\sf M_{41}} & {\sf M_{42}}
\end{matrix}
\right)^{*}
=
\left(
\begin{matrix}
\overline{{\sf M_{31}}} & \overline{{\sf M_{41}}}\\
\overline{{\sf M_{32}}} & \overline{{\sf M_{42}}}
\end{matrix}
\right)
=
-
\left(
\begin{matrix}
{\sf M_{31}} & {\sf M_{32}}\\
{\sf M_{41}} & {\sf M_{42}}
\end{matrix}
\right).
\] 
\subsubsection{Absorption}
The next step is to proceed the absorption process. The point is to simplify the calculations by introducing new variables $W_\bullet$, $X_{\bullet}$, $Y_{\bullet}$ and $Z_{\bullet}$ to modify the coefficents in the matrix $U^{*}dU$: 
\begin{equation}
\begin{aligned}
{\sf M_{31}} &= {\sf \tilde{M}_{31}} + \hspace{1cm} + W_{2}\alpha^{2} + W_{3}\bar{\alpha}^{1}+W_{4}\bar{\alpha}^{2},\\
{\sf M_{32}} &= {\sf \tilde{M}_{32}} + X_{1}\alpha^{1} + \hspace{1cm} + X_{3}\bar{\alpha}^{1}+X_{4}\bar{\alpha}^{2},\\
{\sf M_{41}} &= {\sf \tilde{M}_{41}} + \hspace{1cm} + Y_{2}\alpha^{2} \hspace{0.2cm}+ Y_{3}\bar{\alpha}^{1}+Y_{4}\bar{\alpha}^{2},\\
{\sf M_{42}} &= {\sf \tilde{M}_{42}} + Z_{1}\alpha^{1} \hspace{0.05cm}+ \hspace{1cm} + Z_{3}\bar{\alpha}^{1}+Z_{4}\bar{\alpha}^{2}.
\end{aligned}
\end{equation}
Substituting the equations above into equation \eqref{eq:general-3}, the $2$-forms $d\tilde{\omega}^{3}$ and $d\tilde{\omega}^{4}$ are therefore accordingly changed:
\begin{equation*}
\begin{aligned}
d\tilde{\omega}^{3} &= 
{\sf \tilde{M}_{31}}\wedge\alpha^{1}+{\sf \tilde{M}_{31}}\wedge\alpha^{2}\\
& \hspace{0.5cm} +
(A-W_{3})\ \alpha^{1}\wedge\bar{\alpha}^{1}
+
(B-W_{4})\ \alpha^{1}\wedge\bar{\alpha}^{2}
+
(C-X_{3})\ \alpha^{2}\wedge\bar{\alpha}^{1}
+
(D-X_{4})\ \alpha^{2}\wedge\bar{\alpha}^{2}\\
& \hspace{0.5cm}+
(E-W_{2}+X_{1})\ \alpha^{1}\wedge\alpha^{2},\\
d\tilde{\omega}^{4} &= 
{\sf \tilde{M}_{41}}\wedge\alpha^{1}
+
{\sf \tilde{M}_{42}}\wedge\alpha^{2}\\
& \hspace{0.5cm}+
(F-Y_{3})\ \alpha^{1}\wedge\bar{\alpha}^{1}
+
(G-Y_{4})\ \alpha^{1}\wedge\bar{\alpha}^{2}
+
(H-Z_{3})\ \alpha^{2}\wedge\bar{\alpha}^{1}
+
(I-Z_{4})\ \alpha^{2}\wedge\bar{\alpha}^{2}\\
& \hspace{0.5cm} +
(J-Y_{2}+Z_{1})\ \alpha^{1}\wedge\alpha^{2}.
\end{aligned}
\end{equation*} 
This leads to a system of linear equations
\begin{equation*}
\begin{aligned}
W_{3} &= A,\\
W_{4} &= B,\\
X_{3} &= C,\\
X_{4} &= D,\\
W_{2}-X_{1} &= E,
\end{aligned}
\qquad 
\begin{aligned}
Y_{3} &= F,\\
Y_{4} &= G,\\
Z_{3} &= H,\\
Z_{4} &= I,\\
Y_{2}-Z_{1} &= J.
\end{aligned}
\end{equation*}
The only equations with multiple solutions are $W_{2}-X_{1} = E$ and $Y_{2}-Z_{1} = J$, yet the solutions exist. In the rest of the paper, let $W_{2}^{F}$, $X_{1}^{F}$, $Y_{2}^{F}$ and $Z_{1}^{F}$ denote such a fixed set of solutions.
The absorption continues to take place to obtain a simplification of the Pfaffian system:
\begin{equation} 
\begin{aligned}
d\tilde{\omega}^{3} &= {\sf \tilde{M}_{31}}\wedge\alpha^{1}+{\sf \tilde{M}_{32}}\wedge\alpha^{2},\\
d\tilde{\omega}^{4} &= {\sf \tilde{M}_{41}}\wedge\alpha^{1}+{\sf \tilde{M}_{42}}\wedge\alpha^{2},
\end{aligned}
\end{equation}
with $\varphi^{*}d\tilde{\omega}^{3}=\varphi^{*}d\tilde{\omega}^{4}=0$. 
\subsubsection{Pullback to $\mathbb{D}\times\mathbb{D}$}
The pullback of $d\tilde{\omega}^{3}$ and $d\tilde{\omega}^{4}$ to $\mathbb{D}^{2}$ leads to the vanishing of the differential $2$-forms:
\begin{equation}\label{pullback-eq-1}
\begin{aligned}
0 &= \varphi^{*}{\sf \tilde{M}_{31}}\wedge\delta + \varphi^{*}{\sf \tilde{M}_{32}}\wedge \varepsilon,\\
0 &= \varphi^{*}{\sf \tilde{M}_{41}}\wedge\delta + \varphi^{*}{\sf \tilde{M}_{42}}\wedge\varepsilon.
\end{aligned}
\end{equation}
By Cartan's lemma, there exist functions $k$, $l_{1}$, $l_{2}$ and $m$ such that 
\begin{equation*}
\begin{aligned}
\varphi^{*}{\sf \tilde{M}_{31}} &= k\ \ \delta + l_{1}\ \ \  \varepsilon,\\
\varphi^{*}{\sf \tilde{M}_{32}} &= l_{2}\ \ \delta + m\ \ \varepsilon,
\end{aligned}
\end{equation*}
while in a similar manner, there exist functions $p$, $q_{1}$, $q_{2}$ and $r$ with
\begin{equation*}
\begin{aligned}
\varphi^{*}{\sf \tilde{M}_{41}} &= p\ \ \delta + q_{1}\ \ \varepsilon,\\
\varphi^{*}{\sf \tilde{M}_{42}} &=q_{2}\ \ \delta + r\ \ \varepsilon.
\end{aligned}
\end{equation*}
From equation \eqref{pullback-eq-1}, it follows that $l_{1}=l_{2}$ and $q_{1}=q_{2}$. Let $l:=l_{1}=l_{2}$ and $q:=q_{1}=q_{2}$. Writing
\begin{equation*}
\begin{aligned}
\varphi^{*}A &:= \aaux,\\
\varphi^{*}B &:= \baux,\\
\varphi^{*}C &:= \caux,\\
\varphi^{*}D &:= \daux,
\end{aligned}
\qquad
\begin{aligned}
\varphi^{*}E &:= \eaux,\\
\varphi^{*}F &:= \faux,\\
\varphi^{*}G &:= \gaux,\\
\varphi^{*}H &:= \haux,
\end{aligned}
\qquad
\begin{aligned}
\varphi^{*}I &:= \iaux,\\
\varphi^{*}J &:= \jaux,\\
\varphi^{*}W_{i}&:= \waux_{i},\\
\varphi^{*}X_{i} &:= \xaux_{i},
\end{aligned}
\qquad 
\begin{aligned}
\varphi^{*}Y_{i} &:= \yaux_{i},\\
\varphi^{*}Z_{i} &:= \zaux_{i},\\
\varphi^{*}{\sf M_{ij}} &:= {\sf m_{ij}},\\
~
\end{aligned}
\end{equation*}
therefore,
\begin{equation*}
\begin{aligned}
{\sf m_{31}}  - \waux_{2}\ \varepsilon - \aaux\ \bar{\delta} - \baux\ \bar{\varepsilon} &= k\ \delta +l\ \varepsilon,\\
{\sf m_{32}} -\xaux_{1}\ \delta - \caux\ \bar{\delta} - \daux\ \bar{\varepsilon} &= l\ \delta+m\ \varepsilon,\\
{\sf m_{41}} - \yaux_{2}\ \varepsilon - \faux\ \bar{\delta} - \gaux\ \bar{\varepsilon} &= p\ \delta+q\ \varepsilon,\\
{\sf m_{42}} - \zaux_{1}\ \delta - \haux\ \bar{\delta} - \iaux\ \bar{\varepsilon} &= q\ \delta+r\ \varepsilon.
\end{aligned}
\end{equation*}
In other words,
\begin{equation*}
\begin{aligned}
-{\sf m_{31}} &= -k\hspace{1.4cm}\delta + (-\waux_{2}-l)\ \varepsilon - \aaux\ \bar{\delta} - \baux\ \bar{\varepsilon},\\
-{\sf m_{32}} &= (-\xaux_{1}-l)\hspace{0.25cm} \delta -m\hspace{1.45cm}\ \varepsilon - \caux\ \bar{\delta} - \daux\ \bar{\varepsilon},\\
-{\sf m_{41}} &= -p\hspace{1.3cm}\ \delta + (-\yaux_{2}-q)\ \varepsilon - \faux\ \ \bar{\delta} - \gaux\ \bar{\varepsilon},\\
-{\sf m_{42}} &= (-\zaux_{1}-q)\hspace{0.2cm} \delta -r\hspace{1.6cm}\ \varepsilon - \haux\ \bar{\delta} - \iaux\ \bar{\varepsilon}.
\end{aligned}
\end{equation*}
By the skew-hermitian conditions:
\[
\overline{-{\sf m_{31}}}={\sf m_{31}},
\qquad
\overline{-{\sf m_{32}}}={\sf m_{41}},
\qquad 
\overline{-{\sf m_{41}}} = {\sf m_{32}},
\qquad 
\overline{-{\sf m_{42}}} = {\sf m_{42}},
\]
the following additional restrictions on the torsions are therefore obtained:
\begin{equation*}
\begin{aligned}
& -\overline{k}\bar{\delta}+(-\bar{\waux}_{2}-\bar{l})\bar{\varepsilon}-\bar{\aaux}\delta - \bar{\baux}\varepsilon\\
&= \overline{-{\sf m_{31}}}\\
&= {\sf m_{31}}\\
&= 
k\delta + (\waux_{2}+l)\varepsilon +\aaux\bar{\delta}+\baux\bar{\varepsilon},
\end{aligned}
\qquad
\begin{aligned}
& (-\bar{\xaux}_{1}-\bar{l})\bar{\delta}-\bar{m}\bar{\varepsilon}-\bar{\caux}\delta-\bar{\daux}\varepsilon\\
&= \overline{-{\sf m_{32}}}\\
&= {\sf m_{41}}\\
&=
p\delta+(\yaux_{2}+q)\varepsilon+\faux\bar{\delta}+\gaux\bar{\varepsilon},
\end{aligned}
\qquad
\begin{aligned}
& (-\bar{\zaux}_{1}-\bar{q})\bar{\delta}-\bar{r}\bar{\varepsilon}-\bar{\haux}\delta-\bar{\iaux}\varepsilon\\
&= \overline{-{\sf m_{42}}}\\
&= {\sf m_{42}}\\
&= (\zaux_{1}+q)\delta+r \varepsilon+\haux\bar{\delta}+\iaux\bar{\varepsilon}.
\end{aligned}
\end{equation*}
Finally an inspection of the coefficients yields the set of values:

\begin{equation*}
\begin{aligned}
k &=- \bar{\aaux},\\
l &= -\bar{\baux}-\waux_{2},\\
l &= -\bar{\faux}-\xaux_{1},\\
m &= -\bar{\gaux},\\
p &=- \bar{\caux},\\
q&=- \bar{\daux}-\yaux_{2},\\
q&= -\bar{\haux}-\zaux_{1},\\
r &=- \bar{\iaux}.
\end{aligned}
\end{equation*}
This leads to the question of compatibility amongst the $l$ and $q$:
\[
\bar{\faux}+\xaux_{1} = \bar{\baux}+\waux_{2},\qquad \bar{\haux}+\zaux_{1}=\bar{\daux}+\yaux_{2},
\]
where it is recalled that $\waux_{2}-\xaux_{1}=\eaux$ and $\yaux_{2}-\zaux_{1}=\jaux$. For the holomorphic immersion of the bi-disk to $M^{9}$ passing through the origin to take place, a necessary condition would therefore be
\begin{equation}
\begin{aligned}
\bar{\baux}+\eaux-\bar{\faux}=0\qquad
\text{and}\qquad
\bar{\daux}+\jaux-\bar{\haux}=0.
\end{aligned}
\end{equation}
The following theorem summarises the discussion above:

\begin{Theorem}
	Let $M^{9}\subset\mathbb{C}^{5}$ be a real-analytic smooth real hypersurface passing through the origin, and whose Levi form has signature of $(2,2)$ at each point in $M^{9}$. Suppose $\varphi:\mathbb{D}^{2}\rightarrow M^{9}$ is a holomorphic immersion of bi-disk into $M^{9}$ such that $\varphi(0)=0$, then its unique lift $\tilde{\varphi}:\mathbb{D}^{2}\rightarrow M^{9}\times U(2)$ has the image lying in the zero set of two complex valued functions ${\sf T_{1}}:=\bar{B}+E-\bar{F}=0$ and ${\sf T_{2}}:=\bar{D}+J-\bar{H}=0$.
	\end{Theorem}


\Section{\bf Equivalence method problem under certain conditions}\label{equiv-method}

In this section, to facilitate and simplify the discussion of the equivalence problem, the following hypothesis on the torsions will be assumed:
\begin{equation}\label{twotorsionscondition}
{\sf T_{1}}=\bar{B}+E-\bar{F}\equiv 0,
\qquad\text{and}\qquad 
{\sf T_{2}}=\bar{D}+J-\bar{H}\equiv 0.
\end{equation}
Adopting the notation from the previous section, let $W_{2}^{F}$, $X_{1}^{F}$, $Y_{2}^{F}$ and $Z_{1}^{F}$ be a set of particular solutions to the linear equations
\[
W_{2}^{F}-X_{1}^{F}=E,
\qquad
Y_{2}^{F}-Z_{1}^{F}=J.
\]
The conditions \eqref{twotorsionscondition} then ensure that the two functions $L$ and $Q$ given by
\[
\bar{B}+W_{2}^{F}=: -L :=\bar{F}+X_{1}^{F},
\qquad
\text{and}
\qquad
\bar{D}+Y_{2}^{F}
=:
-Q
:=
\bar{H}+Z_{1}^{F},
\]
are well defined. 
Let $\tau$ be the 2-by-2 matrix of $1$-forms
\begin{equation*}
\begin{aligned}
\tau &:= 
\left(
\begin{matrix}
\hat{\sf M}_{\sf 31}
&
\hat{\sf M}_{\sf 32}\\
\hat{\sf M}_{\sf 41} &
\hat{\sf M}_{\sf 42}
\end{matrix}
\right)\\
&:=
{\footnotesize \left(
\begin{matrix}
-\bar{\sf P}{\sf dP}- \bar{\sf R}{\sf dR}-
W_{2}^{F}\alpha^{2}-A\bar{\alpha}^{1}-B\bar{\alpha}^{2}
+\bar{A}\alpha^{1}-L\alpha^{2}
&
-\bar{\sf P}{\sf dQ}- \bar{\sf R}{\sf dS}-
X_{1}^{F}\alpha^{1}-C\bar{\alpha}^{1}-D\bar{\alpha}^{2}
-L\alpha^{1}+\bar{G}\alpha^{2}\\
-\bar{\sf Q}{\sf dP}- \bar{\sf S}{\sf dR}-
Y_{2}^{F}\alpha^{2}-F\bar{\alpha}^{1}-G\bar{\alpha}^{2}
+\bar{C}\alpha^{1}-Q\alpha^{2}
&
-\bar{\sf Q}{\sf dQ}- \bar{\sf S}{\sf dS}-
Z_{1}^{F}\alpha^{1}-H\bar{\alpha}^{1}-I\bar{\alpha}^{2}
-Q\alpha^{1}+\bar{I}\alpha^{2}
\end{matrix}
\right)}\\
&= 
\left(
\begin{matrix}
-\bar{\sf P}d{\sf P}
-\bar{\sf R} d{\sf R}
+\bar{A}\alpha^{1}
+\bar{B}\alpha^{2}
-A\bar{\alpha}^{1}
-B\bar{\alpha}^{2} &
-\bar{\sf P}d{\sf Q}
-\bar{\sf R}d{\sf S}
+\bar{F}\alpha^{1}
+\bar{G}\alpha^{2}
-C\bar{\alpha}^{1}
-D\bar{\alpha}^{2}\\
-\bar{\sf Q}d{\sf P}
-\bar{\sf S}d{\sf R}
+\bar{C}\alpha^{1}
+\bar{D}\alpha^{2}
-F\bar{\alpha}^{1}
-G\bar{\alpha}^{2}
&
-\bar{\sf Q}d{\sf Q}
-\bar{\sf S}d{\sf S}
+\bar{H}\alpha^{1}
+\bar{I}\alpha^{2}
-H\bar{\alpha}^{1}
-I\bar{\alpha}^{2}
\end{matrix}
\right).
\end{aligned}
\end{equation*}

This matrix is clearly skew-hermitian:
\[
\tau^{*}=-\tau,
\]
and the conditions on the torsions in equation \eqref{twotorsionscondition} also ensure that $\varphi^{*}\tau\equiv 0$.

Now consider any other holomorphic immersion of the bi-disk $\psi:\mathbb{D}^{2}\rightarrow M^{9}$ which passes through the origin, and which satisfies additional conditions that 
\[
\psi^{*}(\omega^{0})\equiv 
\psi^{*}(\omega^{3})\equiv
\psi^{*}(\omega^{4})\equiv 
\psi^{*}(\bar{\omega}^{3})\equiv 
\psi^{*}(\bar{\omega}^{4})\equiv 0,
\qquad
\text{and}
\qquad 
\psi^{*}(\alpha^{1}\wedge\bar{\alpha}^{1}\wedge\alpha^{2}\wedge\bar{\alpha}^{2})\neq 0
\]
at all points in $M^{9}$. Applying  the exterior differentiation leads to the same vanishing of the matrix of $2$-forms
\[
0
\equiv \left(
\begin{matrix}
\psi^{*}d\tilde{\omega}^{3}\\
\psi^{*}d\tilde{\omega}^{4}
\end{matrix}
\right)
=
\psi^{*}
\left[
\left(
\begin{matrix}
{\sf \tilde{M}_{31}} & {\sf \tilde{M}_{32}}\\
{\sf \tilde{M}_{41}} & {\sf \tilde{M}_{41}}
\end{matrix}
\right)
\wedge
\left(
\begin{matrix}
\alpha^{1}\\
\alpha^{2}
\end{matrix}
\right)
\right]
=\psi^{*}\tau  \wedge
\left(
\begin{matrix}
\psi^{*}\alpha^{1}\\
\psi^{*}\alpha^{2}
\end{matrix}
\right).
\]
In view of Cartan's lemma, this necessitates the introduction of new variables $\hat{\kappa}_{i}$ and $\hat{\sigma}_{i}$ such that 
\begin{equation}\label{I3'}
\begin{aligned}
\omega^{0} &= \theta,\\
\omega^{3} &= \alpha^{3}-{\sf P}\alpha^{1}-{\sf Q}\alpha^{2},\\
\omega^{4} &= \alpha^{4}-{\sf R}\alpha^{1} - {\sf S}\alpha^{2},\\
\hat{\omega}^{5} &= 
-\bar{\sf P}d{\sf P}
-
\bar{\sf R}d{\sf R}
-
W_{2}^{F}\alpha^{2}
-
A\bar{\alpha}^{1}
-
B\bar{\alpha}^{2}
-
\hat{\kappa}_{1}\alpha^{1}
-
\hat{\kappa}_{2}\alpha^{2},\\
\hat{\omega}^{6} &= 
-\bar{\sf P}d{\sf Q}
-
\bar{\sf R}d{\sf S}
-
X_{1}^{F}\alpha^{1}
-
C\bar{\alpha}^{1}
-
D\bar{\alpha}^{2}
-
\hat{\kappa}_{2}\alpha^{1}
-
\hat{\kappa}_{3}\alpha^{2},\\
\hat{\omega}^{7} &= 
-\bar{\sf Q}d{\sf P}
-
\bar{\sf S}d{\sf R}
-
Y_{2}^{F}\alpha^{2}
- 
F\bar{\alpha}^{1}
-
G\bar{\alpha}^{2}
-
\hat{\sigma}_{1}\alpha^{1}
-
\hat{\sigma}_{2}\alpha^{2},\\
\hat{\omega}^{8} 
&=
-\bar{\sf Q}d{\sf Q}
-
\bar{\sf S}d{\sf S}
-
Z_{1}^{F}\alpha^{1}
-
H\bar{\alpha}^{1}
-I\bar{\alpha}^{2}
-\hat{\sigma}_{2}\alpha^{1}
-\hat{\sigma}_{3}\alpha^{2}.
\end{aligned}
\end{equation}
These variables depend on $\psi$ which \`{a} priori we have no control of. Another reason for introducing new variables is to resolve the ambiguity of the fixed solutions to the equations
\[W_{2}-X_{1}=E,\qquad \text{and}\qquad
Y_{2}-Z_{1}=J.
\]
More precisely, given a fixed set of solutions $W_{2}^{F}$, $X_{1}^{F}$, $Y_{2}^{F}$ and $Z_{1}^{F}$, the general solutions to these linear equations are 
\[
\left(
\begin{matrix}
W_{2}\\
X_{1}
\end{matrix}
\right)
=
\left(
\begin{matrix}
W_{2}^{F}+\hat{\kappa}_{2}\\
X_{1}^{F}+\hat{\kappa}_{2}
\end{matrix}
\right)
\qquad
\text{and}
\qquad
\left(
\begin{matrix}
Y_{2}\\
Z_{1}
\end{matrix}
\right)
=
\left(
\begin{matrix}
Y_{2}^{F}+\hat{\sigma}_{2}\\
Z_{1}^{F}+\hat{\sigma}_{2}
\end{matrix}
\right)
\]
where $\hat{\kappa}_{2}$ and $\hat{\sigma}_{2}$ become variables that parametrise all the solutions.

Our goal here is to show that the $1$-form $\tau$ is an integral element that vanishes on any holomorphic bi-disk $\psi(\mathbb{D}^{2})$. This will fix values for $\hat{\kappa}_{i}$ and $\hat{\sigma}_{i}$ and thus no further prolongation is necessary.
\begin{Proposition}
Any holomorphic immersion of the bi-disk $\psi:\mathbb{D}^{2}\rightarrow M^{9}$ passing through the origin with
\[
\psi^{*}\omega^{0}
\equiv 
\psi^{*}\omega^{3}
\equiv 
\psi^{*}\omega^{4}
\equiv 
\psi^{*}\bar{\omega}^{3}
\equiv
\psi^{*}\bar{\omega}^{4}
\equiv 0
\]
and
\[
\psi^{*}(\alpha^{1}\wedge\bar{\alpha}^{1}\wedge\alpha^{2}\wedge\bar{\alpha}^{2})\neq 0
\]
at all points in a neighbourhood of the origin necessarily satisfies $\psi^{*}\tau\equiv 0$.
\end{Proposition}

\begin{proof}
The proof uses the fact that $\tau^{*}=-\tau$. For $i=1,2,3,4$, let $w_{i}$, $x_{i}$, $y_{i}$ and $z_{i}$ be functions on $\mathbb{D}^{2}$ such that 
\[
\psi^{*}\tau
=
\left(
\begin{matrix}
w_{1}\psi^{*}\alpha^{1}+
w_{2}\psi^{*}\alpha^{2}+
w_{3}\psi^{*}\bar{\alpha}^{1}+
w_{4}\psi^{*}\bar{\alpha}^{2} &
x_{1}\psi^{*}\alpha^{1}+
x_{2}\psi^{*}\alpha^{2}+
x_{3}\psi^{*}\bar{\alpha}^{1}+
x_{4}\psi^{*}\bar{\alpha}^{2}\\
y_{1}\psi^{*}\alpha^{1}+
y_{2}\psi^{*}\alpha^{2}+
y_{3}\psi^{*}\bar{\alpha}^{1}+
y_{4}\psi^{*}\bar{\alpha}^{2} &
z_{1}\psi^{*}\alpha^{1}+
z_{2}\psi^{*}\alpha^{2}+
z_{3}\psi^{*}\bar{\alpha}^{1}+
z_{4}\psi^{*}\bar{\alpha}^{2}
\end{matrix}
\right).
\]
By the first hypothesis, 
\begin{equation}\label{vanishing-tau-by-pullback-proof-eq-1}
\begin{aligned}
\left(
\begin{matrix}
0 \\
0
\end{matrix}
\right)
&=
-\psi^{*}\tau \wedge
\left(
\begin{matrix}
\psi^{*}\alpha^{1}\\
\psi^{*}\alpha^{2}
\end{matrix}
\right).
\end{aligned}
\end{equation}
The second hypothesis in the proposition statement implies that the $1$-forms $\psi^{*}\alpha^{1}$, $\psi^{*}\alpha^{2}$, $\psi^{*}\bar{\alpha}^{1}$, $\psi^{*}\bar{\alpha}^{2}$ are linearly independent. Using this observation, from the expression of $\psi^{*}\tau$ and equation \eqref{vanishing-tau-by-pullback-proof-eq-1}, the following equations are obtained:
\begin{equation}
\begin{aligned}
w_{2} &= -x_{1},\\
w_{3} &= w_{4} = x_{3} = x_{4} = 0,
\end{aligned}
\qquad
\begin{aligned}
y_{2} &= -z_{1},\\
y_{3} &= y_{4} = z_{3} = z_{4} = 0.
\end{aligned}
\end{equation}
Furthermore, the anti-hermitian property $\psi^{*}\tau^{*}=-\psi^{*}\tau$ results in other equations
\begin{equation}
\begin{aligned}
w_{1}\psi^{*}\alpha^{1}+w_{2}\psi^{*}\alpha^{2} &= 
-\bar{w}_{1}\psi^{*}\bar{\alpha}^{1} 
-\bar{w}_{2}\psi^{*}\bar{\alpha}^{2},\\
\bar{x}_{1}\psi^{*}\bar{\alpha}^{1}+\bar{x}_{2}\psi^{*}\bar{\alpha}^{2}
&=
-y_{1}\psi^{*}\alpha^{1}-y_{2}\psi^{*}\alpha^{2},\\
z_{1}\psi^{*}\alpha^{1}+z_{2}\psi^{*}\alpha^{2}
&=
-\bar{z}_{1}\psi^{*}\bar{\alpha}^{1}-\bar{z}_{2}\psi^{*}
\bar{\alpha}^{2},
\end{aligned}
\end{equation}
which lead to vanishing of the other coefficients in $\psi^{*}\tau$ by the same observation. This implies that $\psi^{*}\tau\equiv 0$.
\end{proof}

Then let $\mathcal{I}_{+}$ be the differential ideal
 \begin{equation}
 \begin{aligned}
\langle \omega^{0},\omega^{3},\omega^{4},\hat{\omega}^{5},\hat{\omega}^{6},\hat{\omega}^{7},\hat{\omega}^{8},\overline{{\omega}^{3}},\overline{{\omega}^{4}},\overline{\hat{\omega}^{5}},\overline{\hat{\omega}^{6}},\overline{\hat{\omega}^{7}},\overline{\hat{\omega}^{8}}\rangle.
\end{aligned}
\end{equation}

In the equation $\psi^{*}\tau=0$, each $\hat{\sf M}_{\sf ij}$ will be required to lie in the ideal $\mathcal{I}_{+}$, and therefore, there exist functions $\eta_{1}$,\dots,$\eta_{13}$ such that 
\[
{\sf \hat{M}_{ij}}
=
\eta_{1}\omega^{0}+\eta_{2}\omega^{3}+\eta_{4}\hat{\omega}^{4}+\cdots +
\eta_{7}\hat{\omega}^{8}+\eta_{9}\overline{{\omega}^{3}}+\eta_{10}\overline{\hat{\omega}^{10}}+\cdots+
\eta_{13}\overline{\hat{\omega}^{8}}.
\]
By direct inspection of the coefficients, it follows that
\begin{equation}\label{solutiontoequations}
\left(
\begin{matrix}
\hat{\kappa}_{1}\\
\hat{\sigma}_{1}
\end{matrix}
\right)
=
\left(
\begin{matrix}
-\bar{A}\\
-\bar{C}
\end{matrix}
\right),
\qquad
\left(
\begin{matrix}
\hat{\kappa}_{2}\\
\hat{\sigma}_{2}
\end{matrix}
\right)
=
\left(
\begin{matrix}
L\\
Q
\end{matrix}
\right),
\qquad
\left(
\begin{matrix}
\hat{\kappa}_{3}\\
\hat{\sigma}_{3}
\end{matrix}
\right)
=
\left(
\begin{matrix}
-\bar{G}\\
-\bar{I}
\end{matrix}
\right).
\end{equation}

\subsubsection{The condition $\psi^{*}d\tau=0$} For the second condition, write $d\tau$ as
\begin{equation*}
\begin{aligned}
d\tau &= 
-\left(
\begin{matrix}
d\bar{\sf P} & d\bar{\sf R}\\
d\bar{\sf Q} & d\bar{\sf S}
\end{matrix}
\right)
\wedge
\left(
\begin{matrix}
d{\sf P} & d{\sf Q}\\
d{\sf R} & d{\sf S}
\end{matrix}
\right)\\
& \hspace{0.5cm} +
d
\left(
\begin{matrix}
-W_{2}^{F}\alpha^{2}-A\bar{\alpha}^{1}-B\bar{\alpha}^{2}+\bar{A}\alpha^{1}-L\alpha^{2}
&
-X_{1}^{F}\alpha^{1}-C\bar{\alpha}^{1}-D\bar{\alpha}^{2}-L\alpha^{1}+\bar{G}\alpha^{2}\\
-Y_{2}^{F}\alpha^{2}-F\bar{\alpha}^{1}-G\bar{\alpha}^{2}+\bar{C}\alpha^{1}-Q\alpha^{2} 
&
-Z_{1}^{F}\alpha^{1}-H\bar{\alpha}^{1}-I\bar{\alpha}^{2}-Q\alpha^{1}+\bar{I}\alpha^{2}
\end{matrix}
\right)\\
&=
-\left(
\begin{matrix}
d\bar{\sf P} & d\bar{\sf R}\\
d\bar{\sf Q} & d\bar{\sf S}
\end{matrix}
\right)
\wedge
\left(
\begin{matrix}
d{\sf P} & d{\sf Q}\\
d{\sf R} & d{\sf S}
\end{matrix}
\right)+
d
\underbrace{\left(
\begin{matrix}
\bar{A}\alpha^{1}+\bar{B}\alpha^{2}-A\bar{\alpha}^{1}-B\bar{\alpha}^{2}
&
\bar{F}\alpha^{1}+\bar{G}\alpha^{2}-C\bar{\alpha}^{1}-D\bar{\alpha}^{2}\\
\bar{C}\alpha^{1}+\bar{D}\alpha^{2}-F\bar{\alpha}^{1}-G\bar{\alpha}^{2}
&
\bar{H}\alpha^{1}+\bar{I}\alpha^{2}-H\bar{\alpha}^{1}-I\bar{\alpha}^{2}
\end{matrix}
\right)}_{:=\Sigma}\\
&=
-\left(
\begin{matrix}
d\bar{\sf P} & d\bar{\sf R}\\
d\bar{\sf Q} & d\bar{\sf S}
\end{matrix}
\right)
\left(
\begin{matrix}
{\sf P} & {\sf Q}\\
{\sf R} & {\sf S}
\end{matrix}
\right)
\wedge
\left(
\begin{matrix}
\overline{\sf P} & \overline{\sf R}\\
\overline{\sf Q} & \overline{\sf S}
\end{matrix}
\right)
\left(
\begin{matrix}
d{\sf P} & d{\sf Q}\\
d{\sf R} & d{\sf S}
\end{matrix}
\right)+ d\Sigma\\
&=
\left(
\begin{matrix}
\overline{\sf P} & \overline{\sf R}\\
\overline{\sf Q} & \overline{\sf S}
\end{matrix}
\right)
\left(
\begin{matrix}
d{\sf P} & d{\sf Q}\\
d{\sf R} & d{\sf S}
\end{matrix}
\right)
\wedge
\left(
\begin{matrix}
\overline{\sf P} & \overline{\sf R}\\
\overline{\sf Q} & \overline{\sf S}
\end{matrix}
\right)
\left(
\begin{matrix}
d{\sf P} & d{\sf Q}\\
d{\sf R} & d{\sf S}
\end{matrix}
\right)+ d\Sigma\\
&\equiv 
\Sigma\wedge\Sigma+d\Sigma \ \ \ \mod\ \mathcal{I}_{+}.
\end{aligned}
\end{equation*}
We therefore have the following theorem 
\begin{Theorem}
	Under the hypothesis that the two essential torsions identically vanish on $M^{9}\times U(2)$, the 2-form $d\tau$ satisfies the following structure equation
	\begin{equation}
	\begin{aligned}
	d\tau &\equiv 
	d\Sigma+\Sigma\wedge\Sigma 
	\equiv 
	0 \qquad \mod\ \mathcal{I}+\langle \tau\rangle.
	\end{aligned}
	\end{equation}
	\end{Theorem}


\Section{\bf Relation to Chern-Moser}
\label{CM}
In this section, we will discuss about the $1$-form $\tau$ in relation to the Chern-Moser tensors, under the condition that the two essential tensors vanish, in other words,  $\bar{B}-E-\bar{F}\equiv 0$ and $\bar{D}+J-\bar{H}\equiv 0$. Then the $2$-forms $d\tau$ has another expression in terms of the $S$-invariants. 

From the general structural equation \eqref{Chern-Moser-structural-eq}, the matrix $g_{\alpha\bar{\beta}}$ chosen for our problem is
\[
g_{\alpha\bar{\beta}}=
\left(
\begin{matrix}
1 & 0 & 0 & 0\\
0 & 1 & 0 & 0\\
0 & 0 & -1  & 0\\
0 & 0 & 0 & -1
\end{matrix}
\right),
\]
so that the first line becomes
\[
d\alpha^{0}=\isqrt\big(\alpha^{1}\wedge\bar{\alpha}^{1}+
\alpha^{2}\wedge\bar{\alpha}^{2}
-\alpha^{3}\wedge\bar{\alpha}^{3}
-\alpha^{4}\wedge\bar{\alpha}^{4}
\big)
+\alpha^{0}\wedge\phi,
\]
corresponding to the Levi form written in the adapted co-frame.
Using $\omega^{3}=\alpha^{3}-{\sf P}\alpha^{1}-{\sf Q}\alpha^{2}$ and $\omega^{4}=\alpha^{4}-{\sf R}\alpha^{1}-{\sf S}\alpha^{2}$, an application of the Poincar\'{e} exterior differentiation yields 
\begin{equation}
\begin{aligned}
\left(
\begin{matrix}
d\omega^{3}\\
d\omega^{4}
\end{matrix}
\right)
&\equiv 
\bigg[
-\left(
\begin{matrix}
\phi_{1}^{3}+{\sf P}\phi_{3}^{3}+{\sf R}\phi_{4}^{3} & 
\phi_{2}^{3}+{\sf Q}\phi_{3}^{3}+{\sf S}\phi_{4}^{3}\\
\phi_{1}^{4}+{\sf P}\phi_{3}^{4}+{\sf R}\phi_{4}^{4} & 
\phi_{2}^{4}+{\sf Q}\phi_{3}^{4}+{\sf S}\phi_{4}^{4}
\end{matrix}
\right)
-
\left(
\begin{matrix}
d{\sf P} & d{\sf Q}\\
d{\sf R} & d{\sf S}
\end{matrix}
\right)\\
&\hspace{0.5cm}+
\left(
\begin{matrix}
{\sf P} & {\sf Q}\\
{\sf R} & {\sf S}
\end{matrix}
\right)
\left(
\begin{matrix}
\phi_{1}^{1}+{\sf P}\phi_{3}^{1}+{\sf R}\phi_{4}^{1} & 
\phi_{2}^{1}+{\sf Q}\phi_{3}^{1}+{\sf S}\phi_{4}^{1}\\
\phi_{1}^{2}+{\sf P}\phi_{3}^{2}+{\sf R}\phi_{4}^{2} & 
\phi_{2}^{2}+{\sf Q}\phi_{3}^{2}+{\sf S}\phi_{4}^{2}
\end{matrix}
\right)
\bigg]
\wedge
\left(
\begin{matrix}
\alpha^{1}\\
\alpha^{2}
\end{matrix}
\right)\,\mod\mathcal{J}\\
&:=
\left(
\begin{matrix}
\check{\sf M}_{\sf 31} & \check{\sf M}_{\sf 32}\\
\check{\sf M}_{\sf 41} & \check{\sf M}_{\sf 42}
\end{matrix}\right)\wedge
\left(
\begin{matrix}
\alpha^{1}\\
\alpha^{2}
\end{matrix}
\right)\,\mod\mathcal{J},
\end{aligned}
\end{equation}
with
\begin{equation}
\begin{aligned}
{\sf \check{M}}_{\sf 31} &=
-d{\sf P}
-
\phi_{1}^{3}
-
{\sf P}\phi_{3}^{3}
-
{\sf R}\phi_{4}^{3}
+
{\sf P}(\phi_{1}^{1}+{\sf P}\phi_{3}^{1}+{\sf R}\phi_{4}^{1})
+
{\sf Q}(\phi_{1}^{2}+{\sf P}\phi_{3}^{2}+{\sf R}\phi_{4}^{2}),\\
{\sf \check{M}}_{\sf 32}
&=
-d{\sf Q}
-\phi_{2}^{3}
-{\sf Q}\phi_{3}^{3}
-{\sf S}\phi_{4}^{3}
+{\sf P}(\phi_{2}^{1}+{\sf Q}\phi_{3}^{1}+{\sf S}\phi_{4}^{1})
+{\sf Q}(\phi_{2}^{2}+{\sf Q}\phi_{3}^{2}+{\sf S}\phi_{4}^{2}),\\
{\sf \check{M}}_{\sf 41} &=
-d{\sf R}
-\phi_{1}^{4}
-
{\sf P}\phi_{3}^{4}
-
{\sf R}\phi_{4}^{4}
+
{\sf R}(\phi_{1}^{1}+{\sf P}\phi_{3}^{1}+{\sf R}\phi_{4}^{1})
+
{\sf S}(\phi_{1}^{2}+{\sf P}\phi_{3}^{2}+{\sf R}\phi_{4}^{2}),\\
{\sf \check{M}}_{\sf 42}
&=
-d{\sf S}
-
\phi_{2}^{4}
-
{\sf Q}\phi_{3}^{4}
-
{\sf S}\phi_{4}^{4}
+
{\sf R}(\phi_{2}^{1}+{\sf Q}\phi_{3}^{1}+{\sf S}\phi_{4}^{1})
+
{\sf S}(\phi_{2}^{2}+{\sf Q}\phi_{3}^{2}+{\sf S}\phi_{4}^{2}).
\end{aligned}
\end{equation}
Then $\tau$ is given by 
\[
\tau = 
\left(
\begin{matrix}
{\sf \bar{P}} & {\sf \bar{R}}\\
{\sf \bar{Q}} & {\sf \bar{S}}
\end{matrix}
\right)
\left(
\begin{matrix}
\check{\sf M}_{\sf 31} & \check{\sf M}_{\sf 32}\\
\check{\sf M}_{\sf 41} & \check{\sf M}_{\sf 42}
\end{matrix}\right).
\]
Using relations \eqref{Chern-Moser-phi}, we have the following:
\begin{Proposition}
The matrix of 1-forms $\tau$ is skew-hermitian.
\end{Proposition}
When the two conditions on the torsions are met, the pullback of $\tau$  vanishes on the bi-disk. Therefore, the Pfaffian system may be expanded to include the coefficients ${\sf \check{M}}_{\sf ij}$:
\[
\mathcal{J}_{+} := \mathcal{J}+\langle \check{\sf M}_{\sf 31},\check{\sf M}_{\sf 32},\check{\sf M}_{\sf 41},\check{\sf M}_{\sf 42}\rangle.
\]
To calculate $d\tau$ modulo $\mathcal{J}_{+}$, it suffices to  do so for 
\[
d\left(
\begin{matrix}
\check{\sf M}_{\sf 31} & \check{\sf M}_{\sf 32}\\
\check{\sf M}_{\sf 41} & \check{\sf M}_{\sf 42}
\end{matrix}\right).
\]
After a lengthy calculation, the following expressions are obtained modulo $\mathcal{I}+\langle\tau\rangle$:
\begin{equation}\label{CM-equation-1}
\begin{aligned}
d\check{\sf M}_{\sf 31}
&\equiv
-\Phi_{1}^{3}-
{\sf P}\Phi_{3}^{3}
-
{\sf R}\Phi_{4}^{3}
+
{\sf P}(\Phi_{1}^{1}
+
{\sf P}\Phi_{3}^{1}
+
{\sf R}\Phi_{4}^{1})
+
{\sf Q}(\Phi_{1}^{2}
+
{\sf P}\Phi_{3}^{2}
+
{\sf R}\Phi_{4}^{2}),\\
d{\sf \check{M}}_{\sf 32}
&\equiv
-\Phi_{2}^{3}
-{\sf Q}\Phi_{3}^{3}
-{\sf S}\Phi_{4}^{3}
+{\sf P}(\Phi_{2}^{1}+{\sf Q}\Phi_{3}^{1}+{\sf S}\Phi_{4}^{1})
+{\sf Q}(\Phi_{2}^{2}+{\sf Q}\Phi_{3}^{2}+{\sf S}\Phi_{4}^{2}),\\
d{\sf \check{M}}_{\sf 41} &\equiv
-\Phi_{1}^{4}
-
{\sf P}\Phi_{3}^{4}
-
{\sf R}\Phi_{4}^{4}
+
{\sf R}(\Phi_{1}^{1}+{\sf P}\Phi_{3}^{1}+{\sf R}\Phi_{4}^{1})
+
{\sf S}(\Phi_{1}^{2}+{\sf P}\Phi_{3}^{2}+{\sf R}\Phi_{4}^{2}),\\
d{\sf \check{M}}_{\sf 42}
&\equiv
-
\Phi_{2}^{4}
-
{\sf Q}\Phi_{3}^{4}
-
{\sf S}\Phi_{4}^{4}
+
{\sf R}(\Phi_{2}^{1}+{\sf Q}\Phi_{3}^{1}+{\sf S}\Phi_{4}^{1})
+
{\sf S}(\Phi_{2}^{2}+{\sf Q}\Phi_{3}^{2}+{\sf S}\Phi_{4}^{2}).
\end{aligned}
\end{equation}
Using the formulae for raising and lowering indices,
\begin{equation}\label{CM-equation-2}
\smallsum{\beta}u_{\alpha}^{\,\,\beta}\ g_{\beta\bar{\gamma}}=u_{\alpha\bar{\gamma}},\qquad
\smallsum{\gamma}u^{\beta}_{\,\,\alpha}\ g^{\alpha\bar{\gamma}}=u^{\beta\bar{\gamma}},
\end{equation}
and with the matrix $g_{\alpha\bar{\beta}}$ being diagonal, the following expressions are obtained
\[
\Phi_{\alpha}^{3}=-\Phi_{\alpha \overline{3}},\qquad 
\Phi_{\alpha}^{4}=-\Phi_{\alpha \overline{4}},\qquad
\Phi_{\alpha}^{1}=\Phi_{\alpha \overline{1}},\qquad
\Phi_{\alpha}^{2}=\Phi_{\alpha \overline{2}}.
\]
Using Lemma 4.2 in \cite{Chern-Moser-1974}, 
\[
\Phi_{\alpha\bar{\rho}}
=
\smallsum{1\leqslant \beta,\bar{\sigma}\leqslant 4}\ S_{\alpha\beta\bar{\rho}\bar{\sigma}}\ 
\alpha^{\beta}\wedge\alpha^{\bar{\sigma}}+\cdots.
\]
Modulo the Pfaffian system which replaces $\alpha^{3}$ by ${\sf P}\alpha^{1}+{\sf Q}\alpha^{2}$, and $\alpha^{4}$ by ${\sf R}\alpha^{1}+{\sf S}\alpha^{2}$, the $2$-forms $\Phi_{\alpha\bar{\rho}}$ may be written as 
\begin{equation}\label{CM-equation-3}
\begin{aligned}
\Phi_{\alpha\bar{\rho}} &\equiv 
\alpha^{1}\wedge\bar{\alpha}^{1}\ 
\big(
S_{\alpha1\bar{\rho}\bar{1}}
+\bar{\sf P}S_{\alpha 1 \bar{\rho}\bar{3}}
+\bar{\sf R}S_{\alpha 1 \bar{\rho}\bar{4}}
+{\sf R}S_{\alpha 4 \bar{\rho}\bar{1}}
+{\sf P}S_{\alpha 3 \bar{\rho}\bar{1}}
+|{\sf P}|^{2}S_{\alpha 3\bar{\rho}3}\\
&\hspace{6 cm}
+{\sf P}\bar{\sf R}S_{\alpha 3 \bar{\rho}\bar{4}}
+\bar{\sf P}{\sf R}S_{\alpha 4 \bar{\rho}3}
+|{\sf R}|^{2}S_{\alpha 4 \bar{\rho}4}
\big)\\
&+
\alpha^{1}\wedge\bar{\alpha}^{2}
\big(
S_{\alpha 1 \bar{\rho}\bar{2}}
+\bar{\sf Q}S_{\alpha 1 \bar{\rho}\bar{3}}
+\bar{\sf S}S_{\alpha 1\bar{\rho}\bar{4}}
+{\sf P}S_{\alpha 3 \bar{\rho}\bar{2}}
+{\sf R}S_{\alpha 4\bar{\rho}\bar{2}}
+{\sf P}\bar{\sf Q}S_{\alpha 3 \bar{\rho}\bar{3}}\\
&\hspace{6cm}
+{\sf P}\bar{\sf S}S_{\alpha 3 \bar{\rho}\bar{4}}
+{\sf R}\bar{\sf Q}S_{\alpha 4 \bar{\rho}\bar{3}}
+{\sf R}\bar{\sf S}S_{\alpha 4 \bar{\rho}\bar{4}}
\big)\\
&+
\alpha^{2}\wedge\bar{\alpha}^{1}
\big(
S_{\alpha 2\bar{\rho}\bar{1}}
+\bar{\sf P}S_{\alpha 2\bar{\rho}\bar{3}}
+\bar{\sf R}S_{\alpha 2 \bar{\rho}\bar{4}}
+{\sf Q}S_{\alpha 3 \bar{\rho}\bar{1}}
+{\sf S}S_{\alpha 4 \bar{\rho}\bar{1}}
+{\sf Q}\bar{\sf P}S_{\alpha 3 \bar{\rho}\bar{3}}\\
&\hspace{6cm}
+{\sf Q}\bar{\sf R}S_{\alpha 3 \bar{\rho}\bar{4}}
+{\sf S}\bar{\sf P}S_{\alpha 4 \bar{\rho}\bar{3}}
+{\sf S}\bar{\sf R}S_{\alpha 4 \bar{\rho}\bar{4}}
\big)\\
&+
\alpha^{2}\wedge\bar{\alpha}^{2}
\big(
S_{\alpha 2 \bar{\rho}\bar{2}}
+\bar{\sf Q}S_{\alpha 2 \bar{\rho}\bar{3}}
+\bar{\sf S}S_{\alpha 2 \bar{\rho}\bar{4}}
+{\sf Q}S_{\alpha 3 \bar{\rho}\bar{2}}
+{\sf S}S_{\alpha 4 \bar{\rho}\bar{2}}
+|{\sf Q}|^{2}S_{\alpha 3 \bar{\rho}\bar{3}}\\
&\hspace{6cm}
+{\sf Q}\bar{\sf S}S_{\alpha 3 \bar{\rho}\bar{4}}
+{\sf S}\bar{\sf Q}S_{\alpha 4 \bar{\rho}\bar{3}}
+|{\sf S}|^{2}S_{\alpha 4 \bar{\rho}\bar{4}}
\big)\ \mod\ \mathcal{I}.
\end{aligned}
\end{equation}
In the expression of $d\tau$:
\[
d\tau\equiv 
\left(
\begin{matrix}
{\sf \bar{P}} & {\sf \bar{R}}\\
{\sf \bar{Q}} & {\sf \bar{S}}
\end{matrix}
\right)
\left(
\begin{matrix}
d\check{\sf M}_{\sf 31} & d\check{\sf M}_{\sf 32}\\
d\check{\sf M}_{\sf 41} & d\check{\sf M}_{\sf 42}
\end{matrix}\right)\qquad\mod\ \mathcal{I}+\langle\tau\rangle,
\]
the use of equations \eqref{CM-equation-1}, \eqref{CM-equation-2} and \eqref{CM-equation-3} leads to 
\[
d\tau=\left(
\begin{matrix}
\sum_{i,j=1}^{2}D_{ij}\alpha^{i}\wedge\overline{\alpha}^{j}
&
\sum_{i,j=1}^{2}E_{ij}\alpha^{i}\wedge\overline{\alpha}^{j}\\
\sum_{i,j=1}^{2}F_{ij}\alpha^{i}\wedge\overline{\alpha}^{j}
&
\sum_{i,j=1}^{2}G_{ij}\alpha^{i}\wedge\overline{\alpha}^{j}
\end{matrix}
\right)
\]
where the $D_{ij}$, $E_{ij}$, $F_{ij}$ and $G_{ij}$ are polynomials with variables ${\sf P}$, ${\sf Q}$, ${\sf R}$, ${\sf S}$, $\bar{\sf P}$, $\bar{\sf Q}$, $\bar{\sf R}$, $\bar{\sf S}$, and with coefficients $S_{ij\bar{k}\bar{l}}$. These polynomials are unfortunately not as simple as the one in the Lorentzian case. For example, the first term $D_{11}$ is given by:
\begin{equation*}
\begin{aligned}
D_{11} &=
S_{33\bar{3}\bar{3}}\ {\sf P}^2\bar{\sf P}^2 
+ 2S_{34\bar{3}\bar{3}}\ {\sf P}{\sf R}\bar{\sf P}^2 
+ S_{44\bar{3}\bar{3}}\ {\sf R}^2\bar{\sf P}^2 
+ 2S_{33\bar{3}\bar{4}}\ {\sf P}^2\bar{\sf P}\bar{\sf R} 
+ 4S_{34\bar{3}\bar{4}}\ {\sf P}{\sf R}\bar{\sf P}\bar{\sf R} 
+ 2S_{44\bar{3}\bar{4}}\ {\sf R}^2\bar{\sf P}\bar{\sf R}\\ 
&+ S_{33\bar{4}\bar{4}}\ {\sf P}^2\bar{\sf R}^2 
+ 2S_{34\bar{4}\bar{4}}\ {\sf P}{\sf R}\bar{\sf R}^2 
+ S_{44\bar{4}\bar{4}}\ {\sf R}^2\bar{\sf R}^2 
+ 2S_{33\bar{1}\bar{3}}\ {\sf P}^2\bar{\sf P} 
+ 4S_{34\bar{1}\bar{3}}\ {\sf P}{\sf R}\bar{\sf P} 
+ 2S_{44\bar{1}\bar{3}}\ {\sf R}^2\bar{\sf P} 
+ 2S_{13\bar{3}\bar{3}}\ {\sf P}\bar{\sf P}^2 \\
&+ 2S_{14\bar{3}\bar{3}}\ {\sf R}\bar{\sf P}^2 
+ 2S_{33\bar{1}\bar{4}}\ {\sf P}^2\bar{\sf R} 
+ 4S_{34\bar{1}\bar{4}}\ {\sf P}{\sf R}\bar{\sf R} 
+ 2S_{44\bar{1}\bar{4}}\ {\sf R}^2\bar{\sf R} 
+ 4S_{13\bar{3}\bar{4}}\ {\sf P}\bar{\sf P}\bar{\sf R} 
+ 4S_{14\bar{3}\bar{4}}\ {\sf R}\bar{\sf P}\bar{\sf R} 
+ 2S_{13\bar{4}\bar{4}}\ {\sf P}\bar{\sf R}^2 \\
&+ 2S_{14\bar{4}\bar{4}}\ {\sf R}\bar{\sf R}^2 
+ S_{33\bar{1}\bar{1}}\ {\sf P}^2 
+ 2S_{34\bar{1}\bar{1}}\ {\sf P}{\sf R} 
+ S_{44\bar{1}\bar{1}}\ {\sf R}^2 
+ 4S_{13\bar{1}\bar{3}}\ {\sf P}\bar{\sf P} 
+ 4S_{14\bar{1}\bar{3}}\ {\sf R}\bar{\sf P} 
+ S_{11\bar{3}\bar{3}}\ \bar{\sf P}^2 
+ 4S_{13\bar{1}\bar{4}}\ {\sf P}\bar{\sf R} \\
&+ 4S_{14\bar{1}\bar{4}}\ {\sf R}\bar{\sf R} 
+ 2S_{11\bar{3}\bar{4}}\ \bar{\sf P}\bar{\sf R} 
+ S_{11\bar{4}\bar{4}}\ \bar{\sf R}^2 
+ 2S_{13\bar{1}\bar{1}}\ {\sf P} 
+ 2S_{14\bar{1}\bar{1}}\ {\sf R} 
+ 2S_{11\bar{1}\bar{3}}\ \bar{\sf P} 
+ 2S_{11\bar{1}\bar{4}}\ \bar{\sf R} 
+ S_{11\bar{1}\bar{1}}.
\end{aligned}
\end{equation*}
The rest of the terms may be found in the arXiv pre-print \cite{Foo-Merker-2Lorentz}.


\Section{\bf An Example}
\label{an-example}

Let $(z_{1},z_{2}, z_{3}, z_{4},u+\isqrt v)$ be holomorphic coordinates of $\mathbb{C}^{5}$.  We will show that if $M^{9}$ is defined by
\[
u=|z_{1}|^{2}+|z_{2}|^{2}-|z_{3}|^{2}-|z_{4}|^{2}+|z_{1}|^{2}(z_{1}+\bar{z}_{1}),
\]
then it does not contain any holomorphic bi-disk passing through the origin by showing that one of the two torsions ${\sf T_{1}}$ or ${\sf T_{2}}$ fail to vanish at the origin.

Consider the following defining equation for $M^{9}$
\[
u=|z_{1}|^{2}+|z_{2}|^{2}-|z_{2}|^{2}-|z_{3}|^{2}+G(z_{1},z_{2},\bar{z}_{1},\bar{z}_{2}),
\]
where $G$ vanishes at the origin at sufficiently high orders. Its Levi matrix is therefore given by
\[
\left(
\begin{matrix}
1+G_{z_{1}\bar{z}_{1}} & G_{z_{1}\bar{z}_{2}} & 0 & 0\\
G_{z_{2}\bar{z}_{1}} & 1+G_{z_{2}\bar{z}_{2}} & 0 & 0\\
0 & 0 & -1 & 0\\
0 & 0 & 0 & -1
\end{matrix}
\right).
\]
To diagonalise the Levi matrix, it suffices to diagonalise the minor matrix
\[
M:=\left(
\begin{matrix}
1+G_{z_{1}\bar{z}_{1}} & G_{z_{1}\bar{z}_{2}} \\
G_{z_{2}\bar{z}_{1}} & 1+G_{z_{2}\bar{z}_{2}} 
\end{matrix}
\right).
\]
This matrix is clearly strictly positive definite in a small neighbourhood of the origin and so it has signature of $(2,2)$. 

\subsection{Change of Coframes}
Let $\alpha^{1}$, $\alpha^{2}$, $\alpha^{3}$ and $\alpha^{4}$ be the $T^{1,0*}M$ frames defined by
\begin{equation*}
\begin{aligned}
\alpha^{1} &= P\ dz^{1}+ Q\ dz^{2},\\
\alpha^{2} &= R\ dz^{2},\\
\alpha^{3} &= dz^{3},\\
\alpha^{4} &= dz^{4}.
\end{aligned}
\end{equation*}
where
\[
P= \left(1+G_{z_{1}\bar{z}_{1}}\right)^{1/2}, 
\qquad
Q=\frac{G_{z_{2}\bar{z}_{1}}}{\left(1+G_{z_{1}\bar{z}_{1}}\right)^{1/2}},
\qquad
R= \frac{\det(M)^{1/2}}{\left(1+G_{z_{1}\bar{z}_{1}}\right)^{1/2}}.
\]
so that the Levi form may be written as
\[
d\alpha^{0} = 
\alpha^{1}\wedge\overline{\alpha}^{1}+\alpha^{2}\wedge\overline{\alpha}^{2}
-\alpha^{3}\wedge\overline{\alpha}^{3}-\alpha^{4}\wedge\overline{\alpha}^{4}.
\]
Moreover, the inverse would be
\begin{equation*}
\begin{aligned}
dz^{1} &= \frac{1}{P}\ \alpha^{1} - \frac{Q}{PR}\ \alpha^{2},\\
dz^{2} &= \frac{1}{R}\ \alpha^{2},\\
dz^{3} &= \alpha^{3},\\
dz^{4} &= \alpha^{4}.
\end{aligned}
\end{equation*}

\subsection{Exterior Derivatives}
Let $z:=(z_{1},\dots,z_{4})$ and let $A(z,\bar{z})$ be any differentiable function. Then 
\begin{equation*}
\begin{aligned}
dA &= \sum_{i=1}^{4}\mathcal{A}_{i}(A)\alpha^{i}+\sum_{i=1}^{4}\bar{\mathcal{A}}_{i}(A)\bar{\alpha}^{i},
\end{aligned}
\end{equation*}
with 
\[
\mathcal{A}_{1}=\frac{1}{P}\frac{\partial}{\partial z_{1}},
\qquad
\mathcal{A}_{2} =- \frac{Q}{PR}\frac{\partial}{\partial z_{1}}+\frac{1}{R}\frac{\partial}{\partial z_{2}},
\qquad
\mathcal{A}_{3} =\frac{\partial}{\partial z_{3}},
\qquad
\mathcal{A}_{4} = \frac{\partial}{\partial z_{4}}.
\]

\subsection{Calculations of $d\alpha^{1}$ and $d\alpha^{2}$}

Since $P$, $Q$ and $R$ depend only on $z_{1}$, $z_{2}$ and their conjugates, the calculations show that 
\begin{equation*}
\begin{aligned}
d\alpha^{1} &= 
\left(
\frac{1}{R}\mathcal{A}_{1}(P)-\frac{Q}{PR}\mathcal{A}_{1}(R)-\frac{1}{P}\mathcal{A}_{2}(P)
\right)\alpha^{1}\wedge\alpha^{2}\\
&
+\left(-\frac{1}{P}\bar{\mathcal{A}}_{1}(P)\right)\alpha^{1}\wedge\bar{\alpha}^{1}
+
\left(\frac{Q}{PR}\bar{\mathcal{A}}_{1}(P)-\frac{1}{R}\bar{\mathcal{A}}_{1}(Q)\right)\alpha^{2}\wedge\bar{\alpha}^{1}\\
&
+\left(-\frac{1}{P}\bar{\mathcal{A}}_{2}(P)\right)\alpha^{1}\wedge\bar{\alpha}^{2}
+
\left(
\frac{Q}{PR}\bar{\mathcal{A}}_{2}(P)-\frac{1}{R}\bar{\mathcal{A}}_{2}(Q)
\right)\alpha^{2}\wedge\bar{\alpha}^{2},\\
d\alpha^{2} &= \frac{1}{R}\mathcal{A}_{1}(R)\alpha^{1}\wedge\alpha^{2}
-\frac{1}{R}\bar{\mathcal{A}}_{1}(R)\alpha^{2}\wedge\bar{\alpha}^{1}-\frac{1}{R}\bar{\mathcal{A}}_{2}(R)\alpha^{2}\wedge\bar{\alpha}^{2},\\
d\alpha^{3} &= 0,\\
d\alpha^{4} &= 0.
\end{aligned}
\end{equation*}

\subsection{The Pfaffian setup}

According to the Pfaffian setup
\begin{equation*}
\begin{aligned}
\omega^{0} &= \theta,\\
\omega^{1} &= \alpha^{1},\\
\omega^{2} &= \alpha^{2},\\
\omega^{3} &= \alpha^{3}-\lambda\alpha^{1}-\mu\alpha^{2},\\
\omega^{4} &= \alpha^{4} - \sigma\alpha^{1} - \theta\alpha^{2},
\end{aligned}
\end{equation*}
with
\[
\left(
\begin{matrix}
\lambda & \mu\\
\sigma & \theta
\end{matrix}
\right)\in U(2),
\]
therefore
\begin{equation*}
\begin{aligned}
 \left(
\begin{matrix}
\bar{\lambda} & \bar{\theta}\\
\bar{\mu} & \bar{\sigma}
\end{matrix}
\right)
\left(
\begin{matrix}
d\omega^{3}\\
d\omega^{4}
\end{matrix}
\right)
&= 
-\left(
\begin{matrix}
\bar{\lambda} & \bar{\theta}\\
\bar{\mu} & \bar{\sigma}
\end{matrix}
\right)
\left(
\begin{matrix}
d\lambda & d\mu\\
d\theta & d\sigma
\end{matrix}
\right)\wedge
\left(
\begin{matrix}
\alpha^{1}\\
\alpha^{2}
\end{matrix}
\right)
+
\left(
\begin{matrix}
\bar{\lambda} & \bar{\theta}\\
\bar{\mu} & \bar{\sigma}
\end{matrix}
\right)
\left(
\begin{matrix}
d\alpha^{3}\\
d\alpha^{4}
\end{matrix}
\right)
-
\left(
\begin{matrix}
d\alpha^{1}\\
d\alpha^{2}
\end{matrix}
\right)\\
&=
-\left(
\begin{matrix}
\bar{\lambda} & \bar{\theta}\\
\bar{\mu} & \bar{\sigma}
\end{matrix}
\right)
\left(
\begin{matrix}
d\lambda & d\mu\\
d\theta & d\sigma
\end{matrix}
\right)
\wedge
\left(
\begin{matrix}
\alpha^{1}\\
\alpha^{2}
\end{matrix}
\right)\\
&\hspace{0.5cm}+
\left(
\begin{matrix}
A\alpha^{1}\wedge\bar{\alpha}^{1}+
B\alpha^{1}\wedge\bar{\alpha}^{2}+
C\alpha^{2}\wedge\bar{\alpha}^{1}+
D\alpha^{2}\wedge\bar{\alpha}^{2}+
E\alpha^{1}\wedge\alpha^{2}\\
F\alpha^{1}\wedge\bar{\alpha}^{1}+
G\alpha^{1}\wedge\bar{\alpha}^{2}+
H\alpha^{2}\wedge\bar{\alpha}^{1}+
I\alpha^{2}\wedge\bar{\alpha}^{2}+
J\alpha^{1}\wedge\alpha^{2}
\end{matrix}
\right),
\end{aligned}
\end{equation*}
with
\begin{equation*}
\begin{aligned}
A &= \frac{1}{P}\bar{\mathcal{A}}_{1}(P),\\
B &= \frac{1}{P}\bar{\mathcal{A}}_{2}(P),\\
C &= \frac{-Q}{PR}\bar{\mathcal{A}}_{1}(P)+\frac{1}{R}\bar{\mathcal{A}}_{1}(Q),\\
D &= \frac{-Q}{PR}\bar{\mathcal{A}}_{2}(P)+\frac{1}{R}\bar{\mathcal{A}}_{2}(Q),\\
E &= \frac{-1}{R}\mathcal{A}_{1}(P)+\frac{Q}{PR}\mathcal{A}_{1}(R)+\frac{1}{P}\mathcal{A}_{2}(P),
\end{aligned}
\qquad
\begin{aligned}
F &= 0,\\
G &= 0,\\
H &= \frac{1}{R}\bar{\mathcal{A}}_{1}(R),\\
I &= \frac{1}{R}\bar{\mathcal{A}}_{2}(R),\\
J &= \frac{-1}{R}\mathcal{A}_{1}(R).
\end{aligned}
\end{equation*}
According to the previous section, the pullback of the following torsions
\begin{equation*}
\begin{aligned}
\bar{B}+E-\bar{F} &= 
\frac{-1}{R}\mathcal{A}_{1}(P)+\frac{Q}{PR}\mathcal{A}_{1}(R)+\frac{2}{P}\mathcal{A}_{2}(P),\\
\bar{D}-J+\bar{H} &=
\frac{-\bar{Q}}{PR}\mathcal{A}_{2}(P)+\frac{1}{R}\mathcal{A}_{2}(\bar{Q})+\frac{2}{R}\mathcal{A}_{1}(R)
\end{aligned}
\end{equation*}
should vanish. 
When $G=|z_{1}|^{2}(z_{1}+\bar{z}_{1})$, 
\[
\bar{B}+E-\bar{F}
=
-\frac{2}{(1+2z_{1}+2\bar{z}_{1})}
\]
which does not vanish near the origin, and therefore, the real hypersurface defined by
\[
u=|z_{1}|^{2}+|z_{2}|^{2}-|z_{3}|^{2}-|z_{4}|^{2}+|z_{1}|^{2}(z_{1}+\bar{z}_{1})
\]
does not contain any bi-disk through the origin.


\bigskip
{\scriptsize
{\sc Wei Guo {\sc Foo}. Hua Loo-Keng Center for Mathematical Sciences, Academy of Mathematics and Systems Science, Chinese Academy of Sciences, Beijing, China.\\
Email address:} \texttt{fooweiguo@hotmail.com}}\\[-10pt]

{\scriptsize
{\sc Jo\"el {\sc Merker}. Laboratoire de Math\'{e}matiques d'Orsay, Universit\'{e} Paris-Sud, CNRS, Universit\'{e} Paris-Saclay, 91405 Orsay Cedex, France.\\
Email address:} \texttt{joel.merker@math.u-psud.fr}}\\[-10pt]


\vfill\end{document}